\documentclass[twoside,11pt,leqno]{article}
\usepackage{amssymb}
\usepackage{amsthm}
\usepackage{amsfonts}
\usepackage[tbtags]{amsmath}
\usepackage{doc}
\usepackage{latexsym}
\font\headd=cmr8

\font\Bbb=msbm10 \theoremstyle{change}

\begin{document}
\small{\addtocounter{page}{208}} \thispagestyle{plain}
 \markboth{}{}
\noindent{\scriptsize KYUNGPOOK Math. J. 40(2000), 209-237}
\vspace{0.2in}

\noindent{\large\bf A Geometrical Theory of Jacobi Forms of Higher
Degree}
 \footnote{(Received : September 9, 1999.\,\,Revised : January 19, 2000) \newline
   \indent \,\,\,\,\,This article is an extended version of the paper published under the same title in Proceedings of Symposium on
        Hodge Theory and Algebraic Geometry (edited by Tadao Oda), Sendai, Japan(1996).}

\vspace{0.15in} \noindent{\sc Jae-Hyun Yang}\newline
 {\it Department of Mathematics, Inha University, Incheon 402-751, Republic of Korea. \\
    e-mail} : {\verb|jhyang@inha.ac.kr|}
\vspace{0.2in}\\
\noindent{\footnotesize (2000 Mathematics Subject Classification : 11F55, 11F46, 11F60, 32M10, 32N10, 53C30)}
 \vspace{0.15in}\\
{\footnotesize \indent In this paper, we give a survey of a
geometrical theory of Jacobi forms of higher degree. And we
present some geometric results and discuss some geometric problems
to be investigated in the future.}

\def\la{\lambda}
\def\a{\alpha}
\def\g{\gamma}
\def\pw{\left({{\partial}\over {\partial W}}\right)}
\def\bhg{\BZ^{(h,g)}_{\geq 0}}
\def\BM{{\mathbb M}(h)}
\def\LM{{\mathcal{L}}_{\mathcal{M}}}
\def\lrt{\longrightarrow}
\def\lmt{\longmapsto}
\def\ba{\backslash}
\def\Hom{\text{Hom}}
\def\s{\sigma}
\def\G{\Gamma}
\def\k{\kappa}
\def\BZ{\mathbb Z}
\def\BQ{\mathbb Q}
\def\BR{\mathbb R}
\def\BC{\mathbb C}
\def\BP{\Bbb P}
\def\vp{\varphi}
\def\J{J\in {\Bbb Z}^{(h,g)}_{\geq 0}}
\def\N{N\in {\Bbb Z}^{(h,g)}}
\def\wt{\widetilde}
\def\Box{$\square$}
\def\M{\mathcal M}
\def\ta{\vartheta^{(\mathcal{M})}\left[\begin{pmatrix} A\\ 0\end{pmatrix}\right]
(\Omega\vert W)}
\def\Jt{\left({{\partial}\over {\partial W}}\right)^J\vartheta^{(\mathcal{M})}
\left[\begin{pmatrix} A\\ 0\end{pmatrix}\right](\Omega\vert W)}
\def\Dt{\Delta^J\vartheta^{(\mathcal{M})}\left[\begin{pmatrix} A\\ 0\end{pmatrix}\right]
(\Omega\vert W)}
\def\emw{\text {exp}\left\{-\pi i\sigma({\mathcal{M}}(\xi\Omega\,^t\xi
+2W\,^t\!\xi))\right\} }
\def\mjz{{\tilde {\vartheta}}^{(\mathcal{M})}_J\left[\begin{pmatrix} A\\ 0\end{pmatrix}
\right](\Omega\vert Z,W)}
\def\CJ{{\Bbb C}\left[\cdots,\left({{\partial}\over
{\partial W}}\right)^J
\vartheta^{(\mathcal{M})}\left[\begin{pmatrix} A_{\alpha,{\mathcal{M}}}\\ 0\end{pmatrix}
\right](\Omega\vert W),\cdots\right]}
\def\GJ{G\left(\cdots,\Delta^J\vartheta^{(\mathcal{M})}\left[\begin{pmatrix}
A_{\alpha,{\mathcal{M}}}\\ 0\end{pmatrix}\right](\Omega\vert W),\cdots\right)}
\def\GW{G\left(\cdots,\left({{\partial}\over {\partial W}}\right)^J
\vartheta^{(\mathcal{M})}\left[\begin{pmatrix} A_{\alpha,{\mathcal{M}}}\\ 0\end{pmatrix}\right]
(\Omega\vert W),\cdots\right)}
\def\lb{\lbrace}
\def\lk{\lbrack}
\def\rb{\rbrace}
\def\rk{\rbrack}
\def\s{\sigma}
\def\w{\wedge}
\def\lrt{\longrightarrow}
\def\lmt{\longmapsto}
\def\O{\Omega}
\def\k{\kappa}
\def\ba{\backslash}
\def\ph{\phi}
\vspace{0.2in}

\noindent{\bf 1. Introduction}
\vspace{0.1in}\\
\indent A Jacobi form is an automorphic form on the Jacobi group,
which is the semi-direct product of the symplectic group
$Sp(g,\BR)$ and the Heisenberg group $H_{\BR}^{(g,h)}$\,(\,see
section 2\,). Jacobi forms are very useful because they are
closely related to modular forms of half integral weight and the
theory of the moduli space of abelian varieties. The simplest case
is when the symplectic group is $SL(2,\BR)$ and the Heisenberg
group is three dimensional, that is, $g=h=1.$ This case had been
treated more or less systematically in [21] and many papers of
Zagier's school. But it seems to us that there is no systematic
investigation of Jacobi forms of higher degree when $g>1$ and
$h>1.$ Some results could be found in [17], [79]-[89] and [94].
\vspace{0.1in}\\
\indent The purpose of this paper is to give a survey of
a geometrical theory of Jacobi forms of higher degree. And we
present some geometric results and discuss some geometric problems
which should be investigated in the future. In Section 2, we
review the notion of Jacobi forms and establish the notations. In
Section 3, we present a brief historical remark and some
motivation on Jacobi forms. In Section 4, we review the toroidal
compactifications of the Siegel modular variety and the universal
abelian variety. In Section 5, we introduce the automorphic vector
bundle $E_{\rho,\M}$ associated with the canonical automorphic
factor $J_{\M,\rho}$ for the Jacobi group $G_{g,h}^J$ and then
discuss the properties of $E_{\rho,\M}$ related to Jacobi forms.
In Section 6, we give some open problems related to Wang's
result(cf. [63]). In Section 7, we describe the boundary of
the Satake compactification in terms of the languages of Jacobi
forms. These results are essentially due to Igusa [35]. In Section
8, we provide you with some characterizations of {\it singular}
Jacobi forms due to Yang [85]. We roughly explain  that the study
of {\it singular} Jacobi forms is closely related to the invariant
theory of the action of the group $GL(g,\BR)\ltimes
H_{\BR}^{(g,h)}$(cf. (9.1)) and to the geometry of the
universal abelian variety. In Section 9, we introduce some results
of the Siegel-Jacobi operator. We describe implicitly that the
Siegel-Jacobi operator plays an important role in the study of the
universal abelian variety. In Section 10, we present
$G_{g,h}^J$-invariant K{\"a}hler metrics and $G_{g,h}^J$-invariant
differential operators on the Siegel-Jacobi space $H_g \times
\BC^{(h,g)}$. We introduce the notion of Maass-Jacobi forms. In
the final section, we give a brief remark on some recent geometric
results. In appendix A, we talk about subvarieties of the Siegel
modular variety and present several problems. In appendix B, we
describe why the study of {\it singular\ modular\ forms} is
closely related to that of the geometry of the Siegel modular
variety. Finally I would like to give my hearty thanks to
Professor Tadao Oda and Dr. Hiroyuki Ito for inviting me to Sendai
and giving me a chance to give a lecture at the conference on
Hodge Theory and Algebraic Geometry.
\vspace{0.1in}\\
\noindent{\bf Notations:}\quad We denote by $\BZ,\,\BR$ and $\BC$
the ring of integers, the field of real numbers, and the field of
complex numbers respectively. $H_g$ denotes the Siegel upper half
plane of degree $g.$ $\G_g:=Sp(g,\BZ)$ denotes the Siegel modular
group of degree $g$. The symbol ``:='' means that the expression
on the right is the definition of that on the left. We denote by
${\BZ}^+$ the set of all positive integers. $F^{(k,l)}$ denotes
the set of all $k\times l$ matrices with entries in a commutative
ring $F$. For a square matrix $A\in F^{(k,k)}$ of degree $k$,
$\sigma(A)$ denotes the trace of $A$. For $A\in F^{(k,l)}$ and
$B\in F^{(k,k)},$ we set $B[A]=\,^t\!ABA.$ For any $M\in
F^{(k,l)},\ ^t\!M$ denotes the transpose matrix of $M$. $E_n$
denotes the identity matrix of degree $n$. For a commutative ring
$K$, we denote by $S_{\ell}(K)$ the vector space of symmetric
matrices of degree $\ell$ with entries in $K$. For a positve
integer $g$ and an integer $k$, we denote by $[\G_g,k]$ the vector
space of all Siegel modular forms on $H_g$ of weight $k$.
 \pagestyle{myheadings}
 \markboth{\headd  Jae-Hyun Yang$~~~~~~~~~~~~~~~~~~~~~~~~~~~~~~~~~~~~~~~~~~~~~~~~~~~~~~~~~~~$}
{\headd $~~~~~~~~~~~~~~~~~~~~~~~~~~~$A Geometrical Theory of Jacobi Forms of Higher Degree}
\vspace{0.1in}\\
\noindent{\bf 2. Jacobi Forms}
\setcounter{equation}{0}
\renewcommand{\theequation}{2.\arabic{equation}}
\vspace{0.1in}\\
\indent In this section, we establish the notations and define the concept of
Jacobi forms.
\vspace{0.05in}\\
\indent Let $$Sp(g,\BR)=\lb M\in \BR^{(2g,2g)}\ \vert \ ^t\!MJ_gM= J_g\ \rb$$
be the symplectic group of degree $g$, where
$$J_g:=\left(\begin{pmatrix} 0 & E_g \\
                   -E_g & 0
                   \end{pmatrix}\right).$$
It is easy to see that $Sp(g,\BR)$ acts on $H_g$ transitively
by
$$M<Z>:=(AZ+B)(CZ+D)^{-1},$$
where $M=\left(\begin{pmatrix} A&B\\ C&D\end{pmatrix}\right)\in Sp(g,\BR)$
and $Z\in H_g.$\\
\indent For two positive integers $g$ and $h$, we recall that
the Jacobi group $G^J_{g,h}:=Sp(g,\BR)\ltimes H_{\BR}^{(g,h)}$ is
the semidirect product of the symplectic group $Sp(g,\BR)$ and
the Heisenberg group $H_{\BR}^{(g,h)}$
endowed with the following multiplication law
$$
(M,(\lambda,\mu,\kappa))\cdot(M',(\lambda',\mu',\kappa'))
:=\, (MM',(\tilde{\lambda}+\lambda',\tilde{\mu}+
\mu', \kappa+\kappa'+\tilde{\lambda}\,^t\!\mu'
-\tilde{\mu}\,^t\!\lambda'))$$
with $M,M'\in Sp(g,\BR), (\lambda,\mu,\kappa),\,(\lambda',\mu',\kappa')
\in H_{\BR}^{(g,h)}$
and $(\tilde{\lambda},\tilde{\mu}):=(\lambda,\mu)M'$.
It is easy to see that $G_{g,h}^J$
acts on $H_{g,h}:=H_g\times \BC^{(h,g)}$ transitively by
\begin{equation}
(M,(\lambda,\mu,\kappa))\cdot (Z,W):=(M<Z>,(W+\lambda Z+\mu)
(CZ+D)^{-1}),
\end{equation}
where $M=\left(\begin{pmatrix} A&B\\ C&D\end{pmatrix}\right)
\in Sp(g,\BR),\ (\lambda,\mu,
\kappa)\in H_{\BR}^{(g,h)}$ and $(Z,W)\in H_{g,h}.$
\vspace{0.05in}\\
\indent Let $\rho$ be a rational representation of $GL(g,\BC)$
on a finite dimensional
complex vector space $V_{\rho}.$ Let ${\mathcal M}\in \BR^{(h,h)}$ be a symmetric
half-integral semi-positive definite matrix of degree $h$.
Let $C^{\infty}(H_{g,h},V_{\rho})$ be the algebra of all
$C^{\infty}$ functions on $H_{g,h}$
with values in $V_{\rho}.$ For $f\in C^{\infty}(H_{g,h},
V_{\rho}),$
we define
$$\begin{array}{ll}
  & (f|_{\rho,{\mathcal{M}}}[(M,(\lambda,\mu,\kappa))])(Z,W)\\
(2.2)\quad\quad:=
\,& e^{-2\pi i\sigma({\mathcal{M}}[W+\lambda Z+\mu](CZ+D)^{-1}C)}
\times
e^{2\pi i\sigma({\mathcal{M}}(\lambda Z^t\!\lambda+2\lambda^t\!W+(\kappa+
\mu^t\!\lambda)))} \\
&\times\rho(CZ+D)^{-1}f(M<Z>,(W+\lambda Z+\mu)(CZ+D)^{-1}),
\hskip 3.8cm
\end{array}$$
where $M=\left(\begin{pmatrix} A&B\\ C&D\end{pmatrix}\right)\in Sp(g,\BR),\
(\lambda,\mu,\kappa)\in H_{\BR}^{(g,h)}$
and $(Z,W)\in H_{g,h}.$
\vspace{0.1in}\\
\def\l{\lambda}
\noindent{\bf Definition\ 2.1.}\quad Let $\rho$ and $\M$ be as above. Let
$$H_{\BZ}^{(g,h)}:=\lb \, (\l,\mu,\kappa)\in H_{\BR}^{(g,h)}\, \vert
\, \l,\mu\in \BZ^{(h,g)},\ \kappa\in \BZ^{(h,h)}\, \rb.$$
Let $\G$ be a discrete subgroup of $\G_g$ of finite index.
A {\it Jacobi\ form} of index $\M$ with respect to $\rho$
on $\G$ is a holomorphic
function $f\in C^{\infty}(H_{g,h},V_{\rho})$ satisfying the
following conditions (A) and (B):
\vspace{0.1in}\\
\indent (A) $f|_{\rho,{\mathcal{M}}}[\tilde{\gamma}] = f$ for all $\tilde{\gamma}\in
\Gamma^J := \G \ltimes H_{\BZ}^{(g,h)}$.\vspace{0.1in}\\
\indent (B) $f$ has a Fourier expansion of the following form :
$$f(Z,W) = \sum\limits_{T\ge0\atop \text {half-integral}}
\sum\limits_{R\in \BZ^{(g,h)}}
c(T,R)\cdot e^{{{2\pi i}\over {\l_{\G}}}\,\sigma(TZ)}\cdot
e^{2\pi i\sigma(RW)}$$
with some nonzero integer $\l_{\G}\in \BZ$ and
$c(T,R)\ne 0$ only if $\left(\begin{pmatrix}
{1\over {\l_{\G}}}T & \frac 12R\\
\frac 12\,^t\!R&{\mathcal{M}}\end{pmatrix}\right)
\ge 0$.\\
\indent If $g\geq 2,$ the condition (B) is superfluous by
K{\" o}cher principle\,(\,cf.\,[94] Lemma 1.6). We denote by $J_{\rho,\M}(\G)$
the vector space of all Jacobi forms of index $\M$ with respect to $\rho$ on $\G$.
Ziegler(cf. [94] Theorem 1.8 or [21] Theorem 1.1) proves that the vector space $J_{\rho,\M}(\G)$ is finite dimensional.
For more results on Jacobi forms with $g>1$ and $h>1$, we refer to
[17], [79]-[89] and [94].
\vspace{0.1in}\\
\noindent{\bf Definition 2.2.}\quad A Jacobi form $f\in J_{\rho,\M}(\G)$ is said to
be a {\it cusp} (or {\it cuspidal}) form if $\begin{pmatrix} {1\over {\l_{\G}}}T & \frac 12R\\
\frac 12\,^t\!R & {\mathcal{M}}\end{pmatrix} > 0$ for any $T,\ R$ with $c(T,R)\ne 0.$
A Jacobi form $f\in J_{\rho,\M}(\Gamma)$ is said to be {\it singular} if it admits a Fourier expansion such that
a Fourier coefficient $c(T,R)$ vanishes unless $\text{det}\,\begin{pmatrix} {1\over {\l_{\G}}}T &\frac 12R\\ \frac 12\,^t\!R&
{\mathcal{M}} \end{pmatrix}=0.$
\vspace{0.1in}\\
\noindent{\bf Example 2.3.}\quad Let $S\in \BZ^{(2k,2k)}$ be a symmetric, positive
definite, unimodular even integral matrix and $c\in \BZ^{(2k,h)}.$
We define the theta series
\begin{equation}
\vartheta_{S,c}^{(g)}(Z,W):=\sum_{\l\in \BZ^{(2k,g)}}
e^{\pi i\{\s(S\l Z\,^t\!\l)+2\s(\,^t\!cS\l\,^t\!W)\} },\ \ Z\in H_g,
\ W\in \BC^{(h,g)}.
\end{equation}
We put $\M:=\frac 12 ^t\!cSc.$ We assume that $2k<g+rank\,(\M).$
Then it is easy to see that $\vartheta_{S,c}^{(g)}$ is a singular
Jacobi form in $J_{k,\M}(\G_g)$(cf. [94] p.212).
\vspace{0.1in}\\
\noindent{\bf 3. Historical Remarks}
\setcounter{equation}{0}
\renewcommand{\theequation}{3.\arabic{equation}}
\vspace{0.1in}\\
\indent In this section, we will make  brief historical remarks on Jacobi forms.
\vspace{0.1in}\\
\indent In 1985, the names Jacobi group and Jacobi forms got kind of standard by the
classic book [21] by {\sc Eichler} and {\sc Zagier} to remind of
Jacobi's ``Fundamenta nova theoriae functionum ellipticorum'', which
appeared in 1829 ([36]). Before [21] these objects appeared more or less
explicitly and under different names in the work of many authors.
\vspace{0.1in}\\
\indent In 1969 Pyatetski-Shapiro [52] discussed the Fourier-Jacobi expansion of
Siegel modular forms and the field of modular abelian functions. He gave
the dimension of this field in the higher degree.\vspace{0.1in}\\
About the same time Satake [55]-[56] introduced the notion of ``groups of
Harish-Chandra type'' which are non reductive but still behave well enough
so that he could determine their canonical automorphic factors and
kernel functions.
\vspace{0.1in}\\
\indent Shimura [57]-[58] gave a new foundation of the theory of complex multiplication
of abelian functions using Jacobi theta functions.
\vspace{0.1in}\\
\indent Kuznetsov [41] constructed functions which are almost Jacobi
forms from ordinary elliptic modular functions.
\vspace{0.1in}\\
\indent Starting 1981, Berndt [4]-[6] published some papers which studied the
field of arithmetic Jacobi functions, ending up with a proof of Shimura
reciprocity law for the field of these functions with arbitrary level.
Furthermore he investigated the discrete series for the Jacobi group
$G^J_{g,h}$
and developed the spectral theory for $L^2(\G^J\backslash G^J_{g,h}$) in the
case $g=h=1$([9],[11]). Recently he [10] studied the
$L$-functions and the Whittaker models for the Jacobi forms.
\vspace{0.1in}\\
\indent The connection of Jacobi forms to modular forms was given by Maass,
Andrianov, Kohnen, Shimura, Eichler and Zagier. This connection is
pictured as follows. For $k$ even, we have the following isomorphisms
$$[\G_2,k]^M\,\cong\,J_{k,1}(\G_1)\,\cong\,M_{k-\frac12}^+(\G_0^{(1)}(4))\,
\cong\,[\G_1,2k-2].$$
Here $[\G_2,k]^M$ denotes the Maass's Spezialschar,
$M_{k-\frac12}^+(\G_0^{(1)}(4))$ denotes the Kohnen space and
$[\G_1, 2k-2]$ denotes the vector space consisting of elliptic modular
forms of weight $2k-2$. For a precise
detail, we refer to [42]-[44], [1], [21], [37] and [81].
\vspace{0.1in}\\
\indent In 1982 Tai [60] gave asymptotic dimension formulae for certain spaces of
Jacobi forms for arbitrary $g$ and $h=1$ and used these ones to show that
the moduli $A_g$ of principally polarized abelian varieties of dimension
$g$ is {\it of\ general\ type} for $g\geq 9.$
\vspace{0.1in}\\
\indent Feingold and Frenkel [23] essentially discussed Jacobi forms in the context
of Kac-Moody Lie algebras generalizing the Maass correspondence to higher
level. Gritsenko [30] studied Fourier-Jacobi expansions and a non-commutative
Hecke ring in connection with the Jacobi group.
\vspace{0.1in}\\
\indent After 1985 the theory of Jacobi forms for $g=h=1$ had been
studied more or less systematically by the
Zagier school. A large part of the theory of Jacobi forms of higher degree
was investigated by Dulinski [17], Kramer [40], Yamazaki [69],
Yang [79]-[89] and Ziegler [94].
\vspace{0.1in}\\
\indent There were several attempts to establish $L$-functions in the
context of the Jacobi group by Murase [47]-[48] and Sugano [50]
using the so-called ``Whittaker-Shintani functions''.
\vspace{0.1in}\\
\indent Recently Kramer [40] developed an arithmetic theory of Jacobi forms of higher degree. Runge [54]
discussed some part of the geometry of Jacobi forms for arbitrary
$g$ and $h=1.$ Quite recently T. Arakawa and B. Heim [2] studied
the iterated Petersson scalar product of a diagonal-restricted
real analytic Jacobi Eisenstein series of degree (3,1) against
elliptic Jacobi forms generalizing Garrett's result in the case of
Siegel Eisenstein series of degree 3.
\vspace{0.1in}\\
\indent For a good survey on some motivation and
background for the study of Jacobi forms, we refer to [10].
\vspace{0.1in}\\
\noindent{\bf 4. Review on Toroidal Compactifications of
the Siegel Space and the Universal Abelian Variety}
\setcounter{equation}{0}
\renewcommand{\theequation}{4.\arabic{equation}}
\vspace{0.1in}\\
\indent In this section, we will make a brief review on toroidal
compactification of the Siegel space and the universal abelian variety.
We refer to [3], [22] and [51]
for more detail.
\vspace{0.1in}\\
\noindent{\bf I. A toroidal compactification of the Siegel modular variety}
\vspace{0.1in}\\
\indent First we realize $H_g$ as a bounded symmetric domain
$D_g:=\,\{\,W\in \BC^{(g,g)}\,\vert\ W=\,^tW,\ \ E_g-Z{\bar Z} > 0\,\}$
\,(\,called the generalized unit disc of degree $g$\,) in $S_g(\BC)$ via
the transformation $\Phi:H_g\lrt D_g$ given by
$$\Phi(Z):=(Z-iE_g)(Z+iE_g)^{-1},\ \ \ Z\in H_g.$$
Indeed, it is a Harish-Chandra realization of a homogeneous space. The
inverse $\Phi^{-1}$ of $\Phi$ given by
$$\Phi^{-1}(Z):=\,i\,(E_g+W)(E_g-W)^{-1},\ \ \ W\in D_g$$
is called the {\it generalized\ Cayley\ transformation.}
\vspace{0.1in}\\
\indent Let ${\bar D}_g$ be the topological closure of $D_g$ in $S_g(\BC).$ Then
${\bar D}_g$ is the disjoint union of all boundary components of $D_g.$
Let
$$F_r:=\left\{ \begin{pmatrix} Z_1& 0 \\
0& E_{g-r}\end{pmatrix} \in {\bar D}_g\, \vert  \ Z_1\in D_r\ \right\},
\ \ \ 0\leq r\leq g$$
be the standard rational boundary components of $D_g.$ Then any boundary
component $F$ of $D_g$ is of the form $F=g\cdot F_r$ for some $g\in
Sp(g,\BR)$ and some $r$ with $0\leq r\leq g.$ In addition, if $F$ is a
rational boundary component of $D_g$, then it is of the form
$F=\gamma\cdot F_r$ for some $\gamma\in Sp(g,\BZ)$ and some $r$ with
$0\leq r\leq g.$ We note that $F_0=\left\{\,E_g\,\right\}$ and
$F_g=D_g.$ We set
\begin{equation}
D_g^*:=\cup_{0\leq r\leq g}\,Sp(g,\BZ)\cdot F_r\,.
\end{equation}
Then $D_g^*$ is clearly the union of all rational boundary components of
$D_g$ and is called the {\it rational\ closure} of $D_g.$ We let
$\G_g:=Sp(g,\BZ)$ for brevity. Then we obtain the so-called
Satake-Baily-Borel compactification $A_g^*:=\G_g\backslash D_g^*$ of
$A_g:=\G_g\backslash D_g.$ Let $F$ be a {\it rational} boundary component
of $D_g.$ We denote by $P(F),\,W(F),\,U(F)$ the parabolic subgroup
associated with $F$, the unipotent radical of $P(F)$ and the center of
$W(F)$ respectively. We set $V(F):=W(F)/U(F)$. Since $P(g\cdot F)=
gP(F)g^{-1}$ for $g\in Sp(g,\BR),$ it is enough to investigate the
structures of these groups for the standard rational boundary components
$F_r\,(\,0\leq r\leq g\,).$
\vspace{0.1in}\\
\indent Now we take $F=F_r$ for some $r$ with $0\leq r\leq g.$ We define
$D(F):=U(F)_{\BC}\cdot D_g\subset {\hat D}_g.$ Here ${\hat D}_g:=
B\backslash Sp(g,\BR)_{\BC}$ is the compact dual of $D_g$ with $B$ a
parabolic subgroup of $Sp(g,\BR)_{\BC}.$ It is obvious that
$U(F)_{\BC}\cong S_{g-r}(\BC)$ and $D(F)\cong F\times V(F)\times U(F)_{\BC}$
analytically. We observe that $U(F)$ acts on $D(F)$ as the linear translation
on the factor $U(F)_{\BC}.$ The isomorphism $\varphi: D(F)\lrt
F\times V(F)\times U(F)_{\BC}$ is given by
$$\varphi\left( \begin{pmatrix} Z_1 & Z_2\\
* & Z_3 \end{pmatrix} \right):=\,(Z_1,Z_2,Z_3),\ \ Z_1\in D_r,
\ Z_2\in \BC^{(r,g-r)},\ Z_3\in S_{g-r}(\BC).$$
We define the mapping $\Phi_F: D(F)\lrt U(F)$ by
\begin{equation}
\Phi_F((Z_1,Z_2,Z_3)):=\text{Im}\,Z_3\,-\,^t(\text{Im}\,Z_2)\,
(\text{Im}\,Z)^{-1}\,(\text{Im}\,Z_2),\ \ \ (Z_1,Z_2,Z_3)\in D(F).
\end{equation}
Then $D_g\cong H_g$ is characterized by $\Phi_F(Z) > 0$ for all $Z\in D_g.$
This is the realization of a Siegel domain of the third kind. We let
$C(F)$ be the cone of real positive symmetric matrices of degree $g-r$ in
$U(F)\cong S_{g-r}(\BR).$ Clearly we have $D_g=\Phi^{-1}(C(F)).$ We define
$$G_h(F):=\text{Aut}\,(F)\ \ (\,\text{modulo\ finite\ group}\,)$$
and
$$G_l(F):=\text{Aut}\,(U(F),\,C(F)).$$
Then it is easy to see that
$$P(F)\,=\,\left(\,G_h(F)\times G_l(F)\,\right)\ltimes W(F)\ \ \
(\,\text{the\ semidirect\ product}\,)  $$
We obtain the natural projections $p_h:P(F)\lrt G_h(F)$ and $p_l:P(F)\lrt
G_l(F)$.
\vspace{0.1in}\\
\noindent{\bf Step I} : Partial compactification for a rational bounadry
component.
\vspace{0.1in}\\
\indent Now we let $\G$ be an arithmetic subgroup of $Sp(g,\BR)$. We let
$$\begin{array}{ll}
\G (F):&=\G\cap P(F),\\
{\bar {\G}}(F):&=p_l(\G(F))\subset G_l(F),\\
U_{\G}(F):&=\G\cap U(F),\ \ \ \text{a\ lattice\ in}\ U(F),\\
W_{\G}(F):&=\G\cap W(F).
\end{array}$$
We note that ${\bar {\G}}(F)$ is an arithmetic subgroup of $G_l(F).$
\vspace{0.1in}\\
\indent Let $\Sigma_F=\left\{\,\s_{\a}^F\,\right\}$ be a ${\bar {\G}}(F)$-admissible
polyhedral decomposition of $C(F).$ We set $D(F)':=D(F)/U(F)_{\BC}.$
Since $D(F)'\cong F\times V(F),$ the projection $\pi_F:D(F)\lrt D(F)'$ is
a principal $U(F)_{\BC}$-bundle over $D(F)'.$ The map
\begin{equation}
\pi_{F,\G}:U_{\G}(F)\ba D(F)\cong F\times V(F)\times \left( U_{\G}(F)\ba
U(F)_{\BC}\right)\lrt D(F)'
\end{equation}
is a principal $T(F)$-bundle with the structure group $T(F):=
U_{\G}(F)\ba U(F)_{\BC}\cong (\BC^*)^q,$ where $q=\dfrac{(g-r)(g-r+1)}{2}.$
Let $\text{X}_{\Sigma_F}$ be a normal torus embedding of $T(F).$ We note
that $\text{X}_{\Sigma_F}$ is determined by $\Sigma_F.$ Then we obtain a
fibre bundle
\begin{equation}
{\mathcal X}(\Sigma_F):=\left( U_{\G}(F)\ba D(F)\right)\times_{T(F)}\,
{\text{X}}_{\Sigma_F}
\end{equation}
over $D(F)'$ with fibre $\text{X}_{\Sigma_F}.$ We denote by {\bf X}
$(\Sigma_F)$ the interior of the closure of $U_{\G}(F)\ba D_g$ in
${\mathcal X}(\Sigma_F)$\,(\,because $D_g\subset D(F)$\,). {\bf X}$(\Sigma_F)$
has a fibrewise $T(F)$-orbit decomposition $\coprod_{\mu}O(\mu)$
such that
$$\begin{array}{ll}
\ (i)\ & \text{each}\ O(\mu)\ \text{is\ an\ algebraic\ torus\ bundle\
over}\ D(F)',\\
(ii)\ & \s_{\mu}\prec \s_{\nu}\ \ \text{iff}\ \ {\overline{O(\mu)}}
\supseteq O(\nu),\\
(iii)\ & \text{dim}\,\s_{\mu}\,+\,\text{dim}\,O(\mu)\,=\,\text{dim}\,D(F),\\
(iv)\ & \text{for}\ \s_{\mu}=\,{0},\,\ O(\mu)=U_{\G}(F)\ba D(F).
\end{array}$$
We define
$$O(F):=\bigcup_{\s_{\a}^F\cap C(F)\neq \emptyset}\,O(\a)\subset
{\bold X}(\Sigma_F)$$
and
$${\bar O}(F):=\G(F)/U_{\G}(F)\ba O(F).$$
We note that $O(F_g)=D_g$ and ${\bar O}(F_g)=\G\ba D_g.$
We set
\begin{equation}
{\bold Y}(\Sigma_F):=\G(F)/U_{\G}(F)\ba {\bold X}(\Sigma_F)\,.
\end{equation}
We note that
$\G(F)/U_{\G}(F)$ acts on {\bf Y}$(\Sigma_F)$ properly discontinuously.
Then we can show that {\bf Y}$(\Sigma_F)$ has a canonical quotient structure
of a normal analytic space and ${\bar O}(F)$ is a closed analytic set in
{\bf Y}$(\Sigma_F).$
\vspace{0.1in}\\
\noindent{\bf Step II} : Gluing.
\vspace{0.1in}\\
\indent Let $\Sigma:=\{\,\Sigma_F\,\vert\ F\ \text{is\ a\ rational\ boundary\
component\ of}\ D_g\,\}$ be a $\G$-admissible family of polyhedral
decompositions. We put
$${\widetilde{\left(\G\ba D_g\right)}}:=\cup_{F:\text{rational}} {\bold X}
(\Sigma_F).$$
We define the equivalence relation $\backsim $ on
${\tilde{\left( \G\ba D_g\right)}}$ as follows:
$$X_1\,\backsim\,X_2,\ \ X_1\in {\bold X}(\Sigma_{F_1}),\ \ X_2\in
{\bold X}(\Sigma_{F_2})$$
iff there exist a rational boundary component $F$, an element $\gamma\in
\G$ such that $F_1\prec F,\ \gamma\,F_2\prec F$ and there exists an element
$X\in {\bold X}(\Sigma_F)$ such that $\pi_{F,F_1}(X)=X_1,\ \pi_{F,F_2}(X)=
\gamma X_2,$ where
$$\pi_{F,F_1}:{\bold X}(\Sigma_F)\lrt {\bold X}(\Sigma_{F_1}),\ \ \
\pi_{F,F_2}:{\bold X}(\Sigma_F)\lrt {\bold X}(\Sigma_{\gamma F_2}).$$
The space $\overline{\left( \G\ba D_g\right)}:=
{\widetilde{\left( \G\ba D_g\right)}}/\backsim$ is called the {\it toroidal\
compactification} of $\G\ba D_g$ associated with $\Sigma$. It is known that
$\overline{\left( \G\ba D_g\right)}$ is a Hausdorff analytic variety
containing $\G\ba D_g$ as an open dense subset. For a {\it neat}
arithmetic subgroup $\G$, we can obtain a smooth projective toroidal
compactification of $\G\ba D_g.$
\vspace{0.1in}\\
\noindent{\bf II. A toroidal compactification of the universal abelian variety}
\vspace{0.1in}\\
\indent For a positive integer $g\in \BZ^+,$ we put $X:=\BZ^g.$ Let $B(X)$ be the
$\BZ$-module of integral valued symmetric bilinear forms on $X$ and let
$B(X)_{\BR}:=B(X)\otimes_{\BZ}\BR.$ Let $C(X)\subset B(X)_{\BR}$ be the
convex cone of all positive semi-positive symmetric bilinear forms on
$X_{\BR}$ whose radicals are defined over $\BQ.$ We let $X^{\ast}$ be the
dual of $X$. For a positive integer $s\in \BZ^+,$ we let
$${\tilde B}_s(X):=B(X)\times (X^{\ast})^s\ \ \ \text{and}\ \ \
{\tilde B}_s(X)_{\BR}:={\tilde B}_s(X)\otimes_{\BZ}\BR.$$
Then the semidirect product $GL(X)\ltimes X^s$ acts on
${\tilde B}_s(X)_{\BR}$ in the natural way and the projection
${\tilde B}_s(X)_{\BR}\lrt B(X)_{\BR}$ is equivariant with respect to the
canonical morphism $GL(X)\ltimes X^s\lrt GL(X).$ Inside
${\tilde B}_s(X)_{\BR}$ we obtain the cone ${\tilde C}_s(X)$ consisting of
$q=(b;\ell_1,\cdots,\ell_s)\in {\tilde B}_s(X)_{\BR}$ such that
$b\in C(X)$ and each $\ell_j$ vanishes on the radical of $b$.
\vspace{0.1in}\\
\indent Let a $GL(X)$-admissible polyhedral cone decomposition ${\mathcal{C}}=
\{ \s_{\a}\}$ of $C(X)$ be given. A $GL(X)\ltimes X^s$-admissible
polyhedral cone decomposition
${\tilde {\mathcal{C}}}=\{\tau_{\beta}\}$ of
${\tilde C}_s(X)$ relative to ${\mathcal{C}}=\{\s_{\a}\}$ is defined to be a
collection ${\tilde {\mathcal{C}}}=\{\tau_{\beta}\}$ such that\\
\indent (1) each $\tau_{\beta}$ is a non-degenerate rational polyhedral cone which is
open in the smallest $\BR$-subspace containing it;\\
\indent (2) any face of a $\tau_{\beta}\in {\tilde {\mathcal{C}}}$ belongs to
${\tilde {\mathcal{C}}}$;\\
\indent (3) ${\tilde C}_s(X)=\cup_{\tau_{\beta}\in {\tilde {\mathcal{C}}}}
\tau_{\beta}\,;$\\
\indent (4)\ ${\tilde {\mathcal{C}}}$ is invariant under the action of
$GL(X)\ltimes X^s$ and there are only finitely many
$GL(X)\ltimes X^s$-orbits;\\
\indent (5) any $\tau_{\beta}\in {\tilde {\mathcal{C}}}$ maps into a $\s_{\a}\in
{\mathcal{C}}$ under the natural projection ${\tilde C}_s(X)\lrt C(X).$
\vspace{0.1in}\\
\indent We call ${\tilde {\mathcal{C}}}$ {\it equidimensional} if in (5) of the above
definition each $\tau_{\beta}\in {\tilde {\mathcal{C}}}$ maps onto a
$\s_{\a}\in {\mathcal{C}}.$
\def\C{\mathcal{C}}
\def\TC{\tilde {\mathcal{C}}}
Again, $\TC$ is called {\it smooth} or {\it regular} if each
$\tau_{\beta}\in \TC$ is generated by part of a $\BZ$-basis of
${\tilde B}_s(X).$
\def\BS{{\tilde B}_s(X)}
\def\CS{{\tilde C}_s(X)}
According to the reduction theory [3], there exists
a smooth equidimensional
$GL(X)\ltimes X^s$-admissible polyhedral cone decomposition $\TC$ of
$\CS$ relative to ${\mathcal{C}}.$ Let $F$ be the split torus $\BS_{\BR}
\otimes_{\BZ} {\Bbb G}_{\text{\bf m}}.$ The choice of a polyhedral cone
decomposition $\TC=\{\tau_{\beta}\}$ of $\CS$ as above provides us with
a torus embedding $F\hookrightarrow {\bar F}.$ Then ${\bar F}$ is
stratified by $F$-orbits and $GL(X)\ltimes X^s$ acts on ${\bar F}$
preserving this stratification. Therefore we obtain the toroidal
compactification ${\bar A}_{g,s}$ of the universal abelian variety
$A_{g,s}:=\G_{g,s}^J\ba H_g\times \BC^{(s,g)}$ with
$\G_{g,s}^J:=\G_g\ltimes H_{\BZ}^{(g,h)}.$ We collect some properties of
the toroidal compactification ${\bar A}_{g,s}.$
\vspace{0.1in}\\
\indent (a) ${\bar A}_{g,s}$ is a Haudorff analytic variety containing $A_{g,s}$
as an open dense subset.\\
\indent (b) ${\bar A}_{g,s}$ has a stratification parametrized by the
$GL(X)\ltimes X^s$-orbits of cones $\tau_{\beta}\in \TC.$
\\
\indent (c) The toroidal compactification ${\bar A}_{g,s}$ depends on the choice of
a smooth equidimensional $GL(X)\ltimes X^s$-admissible polyhedral cone
decomposition $\TC=\{\tau_{\beta}\}$ of
$\CS$ relative to ${\mathcal{C}}.$ In order to
indicate this dependence we write ${\bar A}_{g,s}(\TC)$ instead of
${\bar A}_{g,s}.$ The natural projection $\pi:A_{g,s}\lrt A_g$ extends to
a proper morphism ${\bar {\psi}}:{\bar A}_{g,s}\lrt {\bar A}_g.$
\vspace{0.1in}\\
\indent Now we recall [22],\,p.\,197 that an {\it admissible\ homogeneous\
principal\ polarization\ function} of $\{\tau_{\beta}\}\lrt \{\s_{\a}\}$ is
a piecewise linear function ${\tilde{\phi}}:\CS\lrt \BR$ satisfying the
following conditions\vspace{0.1in}\\
\def\TP{\tilde{\phi}}
\indent (P1) $\TP$ is continuous and $GL(X)$-invariant;\\
\indent (P2) $\TP$ takes rational values on $\BS\cap \CS$ with bounded denominators;
\\
\indent (P3) $\TP$ is homogeneous, i.e., $\TP(t\cdot q)=t\cdot \TP(q)$ for all real
$t\geq 0$ and all $q\in \CS$;\\
\indent (P4) $\TP$ is linear on each $\tau_{\beta}\in \TC$;\\
\indent (P5) $\TP$ is convex in the sense that
$$\TP(t\cdot q+(1-t)\cdot q')\geq t\cdot \TP(q)+(1-t)\cdot \TP(q')$$
for all $t\in\BR$ with $0\leq t\leq 1$ and any $q,q'\in \CS.$
\\
\indent (P6) $\TP$ is strictly convex, that is, for each $\s_{\a}\in {\mathcal{C}}=
\{\s_{\a}\}$ and each $\tau_{\beta}\in \TC=\{\tau_{\beta}\}$ lying over
$\s_{\a}$, there exist a finite number of linear functionals
$\ell_i:\BS\lrt \BR,\ 1\leq i\leq m$ with $\ell_i\geq \TP$ on the preimage
of $\s_{\a}$ for each $i$ and
$$\tau_{\beta}=\{\,q\in \CS\,\vert\ q\ \text{lies\ over}\ \s_{\a}\ \text{and}\
\TP(q)=\ell_i(q)\ \text{for\ each}\ i\ \}.$$
\indent (P7) There exists a rational positive number $r$ such that for each
$\mu=(\mu_1,\cdots,\mu_s)\in X^s,$ the function
$$\TP-\TP\circ T_{\mu}\,:\,q\longmapsto f(q)-f(\mu\cdot q)$$
is equal to $r$ times (restriction to $\CS$ of) the linear functional
${\tilde {\chi}}_{\mu}$ on $\BS$, where for $q=(b;\ell_1,\cdots,\ell_s)\in
\CS,$
$${\tilde{\chi}}_{\mu}(q):=\sum_{1\leq i\leq s}a_i(\mu_i)=
\sum_{1\leq i\leq s}\{\,b(\mu_i,\mu_i)+2\cdot \ell_i(\mu_i)\}.$$
The conditions (P1)-(P7) above constitute a kind of convexity conditions
on $\{\tau_{\beta}\}\lrt \{\s_{\a}\}.$ They imply that the morphism
${\bar A}_{g,s}\lrt {\bar A}_g$ attached to $\{\tau_{\beta}\}\lrt
\{\s_{\a}\}$ is {\it projective}. Indeed, the theory of torus embeddings
shows that an admissible homogeneous principal polarization function
$\TP:\CS\lrt \BR$ gives rise to an invertible sheaf
${\bar{\mathcal{L}}}(\TP),$
which is ample
on ${\bar A}_{g,s}(\tilde{\mathcal{C}})$ relative to ${\bar A}_g({\mathcal{C}}).\hfill$
\vspace{0.1in}\\
\noindent{\bf 5. The Automorphic Vector Bundle} $E_{\rho,\M}$
\setcounter{equation}{0}
\renewcommand{\theequation}{5.\arabic{equation}}
\vspace{0.1in}\\
\indent Let $\rho$ and $\M$ be as before in section 2.
Assume that $\Gamma$ is a subgroup of $\G_g:=Sp(g,{\Bbb Z})$ of finte index
which acts freely on $H_g$ and $-E_{2g}\notin \Gamma.$
Then $\G^J:=\G\ltimes
H_{\BZ}^{(g,h)}$ acts on $H_{g,h}:=H_g\times \BC^{(h,g)}$ properly
discontinuously. We consider the automorphic factor
$J_{\M,\rho}:G^J_{g,h}\times
H_{g,h}\lrt GL(V_{\rho})$ defined by
$$\begin{array}{ll}
J_{\M,\rho}({\tilde{g}},(Z,W)):=&e^{2\pi i\s (\M [W+\lambda Z+\mu]
(CZ+D)^{-1}C)}\\
&\times e^{-2\pi i\s(\M(\la Z\,^t\la+2\la ^tW+\kappa+\mu\,^t\la))}\rho(CZ+D),
\end{array}$$
where ${\tilde g}=(M,(\la,\mu,\kappa))\in G^J$ with $M=
\begin{pmatrix} A & B\\ C & D \end{pmatrix} \in Sp(g,\BR).$ Then $J_{\M,\rho}$
defines the automorphic vector bundle $E_{\rho,\M}:=H_{g,h}\times_{\G^J}
V_{\rho}$ over $A_{g,h,\G}:=\G^J\ba H_{g,h}.$ By the definition,
Jacobi forms in $J_{\rho,\M}(\G)$ may be considered as holomorphic sections
of the vector bundle $E_{\rho,\M}$ with some additional cusp condition.
For $g\geq 2,$ this additional condition may be dropped according to
K{\"o}cher principle. Let ${\bar A}_{g,h,\G}$ be a toroidal compactification
given by a regular $\G$-admissible family $\Sigma$ of polyhedral
decompositions.\vspace{0.1in}\\
Without proof we provides our results.
\vspace{0.1in}\\
\noindent{\bf Theorem 5.1.}\quad {\it $A_{g,h,\G}$ is contained in ${\bar A}_{g,h,\G}$
as a Zariski open subset. $E_{\rho,\M}$ can be extended uniquely to
the holomorphic vector bundle ${\bar E}_{\rho,\M}$ over
${\bar A}_{g,h,\M}.$ And $H^i(A_{g,h,\G},E_{\rho,\M})$\\
\noindent
$\cong H^i({\bar A}_{g,h,\M},{\bar E}_{\rho,\M}).$ In particular, the dimension of
$J_{\rho,\M}$ is finite dimensional.}
\\vspace{0.1in}\\
\noindent{\bf Definition 5.2.}\quad Let $\rho$ be an irreducible rational representation
of $GL(g,\BC)$ with its highest weight $(\l_1,\l_2,\cdots,\l_g).$
We call the number of $j\,(\,1\leq j\leq g\,)$ such that $\l_j=\l_g$
the {\it corank} of $\rho$ which is denoted by $\text{corank}\,(\rho).$
The number $k(\rho):=\l_g$ is called the {\it weight} of $\rho.$
\vspace{0.1in}\\
\noindent{\bf Theorem 5.3.}\quad {\it Let $2\M$ be an even unimodular positive definite
matrix of degree $h.$ Let $\rho$ be an irreducible finite dimensional
representation of $GL(g,\BC)$ with highest weight $\rho=(\l_1,\cdots,
\l_g).$ Let $\l(\rho)$ be the number of $\l_i's$ such that
$\l_i=k(\rho)+1=\l_g+1,\ 1\leq i\leq g.$ Assume that $\rho$ satisfies the
following conditions :
$$\begin{array}{ll}
&[a]\quad \rho(A)=\rho(-A)\ \ \ \text{for\ all}\ A\in GL(g,\BC),\\
&[b]\quad \l(\rho) < 2(g-k(\rho)-\text{corank}\,(\rho)\,)+h.
\end{array}$$
Then $H^0(A_{g,h,\G},E_{\rho,\M})=0.$}
\vspace{0.05in}\\
\noindent{\it Proof.}\quad The proof can be found in [80].
\vspace{0.1in}\\
\noindent{\bf Corollary 5.4.}\quad {\it Let $2\M$ be as above in Theorem 5.3. Assume that
$2k(\rho)\leq g+h-2\text{corank}\,(\rho).$ Then $H^0(A_{g,h,\G},
E_{\rho,\M})=0.$}
\vspace{0.1in}\\
\noindent{\bf Remark 5.5.}\quad N.-P. Skoruppa [Sk] proved that $J_{1,m}(\G_1)=0$
for any nonnegative integer $m$. It is interesting to give the geometric
proofs of this fact and Theorem 5.3.
\vspace{0.1in}\\
We give the following open problems :
\vspace{0.1in}\\
\noindent{\bf Problem 1.}\quad Give the explicit dimension formula or estimate for
$H^0(A_{g,h,\G},E_{\rho,\M}).$
\vspace{0.1in}\\
\noindent{\bf Problem 2.}\quad Compute the cohomology groups $H^k(A_{g,h,\G},\,
E_{\rho,\M})$ explicitly. Here $0\leq k\leq \dfrac{g(g+2h+1)}{2}$.
\vspace{0.1in}\\
\noindent{\bf Problem 3.}\quad Under which conditions is $E_{\rho,\M}$ {\it ample}?
\vspace{0.1in}\\
\noindent{\bf Problem 4.}\quad Discuss the analogue of Hirzebruch's proportionality
theorem for $E_{\rho,\M}$(cf. [45]).
\vspace{0.1in}\\
\noindent{\bf 6. Smooth Compactification of Siegel Moduli Spaces and Open Problems}
\setcounter{equation}{0}
\renewcommand{\theequation}{6.\arabic{equation}}
\vspace{0.1in}\\
\indent Let $\G_g(k)$ be the principal congruence subgroup of $Sp(g,{\Bbb Z})$ of
level $k$ and let $H_g$ be the Siegel upper-half plane of degree $g$. We
assume that $k\geq 3.$ This implies that $\G_g(k)$ is a {\it neat} arithmetic
subgroup. Let ${\bar X}$ be the toroidal compactification of $X:=\G_g(k)
\ba H_g$ from $\G_g(k)$-admissible family given by the central cone
decomposition $\sum_{cent}$ or a refinement of $\sum_{cent}$. Then the
boundary $D:={\bar X}-X=\sum_{i=1}^m D_i$ is a divisor of ${\bar X}$ with
normal crossing, that is, each $D_i$ is an irreducible smooth divisor of
${\bar X}$ and $D_1,\cdots,D_m$ intersect transversally. If $g\leq 4,$
we have the following results obtained by Wang [63].
\vskip 0.2cm \noindent
{\bf Theorem 6.1.}\quad {\it (1) Each divisor is algebraically isomorphic to
$${\bar Y}_{g-1}:={\overline {\G_{g-1}^J(k)\ba (H_{g-1}\times
{\Bbb C}^{g-1})}}.$$
Here $\G_{g-1}^J(k):=\G_{g-1}(k)\ltimes (k{\Bbb Z})^{g-1}$
is the Jacobi modular group acting on the homogeneous space
$W_{g-1}:= H_{g-1}\times {\Bbb C}^{g-1}$ in a usual way and ${\bar Y}_{g-1}$
is the compactification of the universal family $Y_{g-1}:=\G_{g-1}^J(k)\ba
(H_{g-1}\times {\Bbb C}^{g-1})$ of abelian varieties induced from the same
$\G_g(k)$-admissible family.
\vspace{0.05in}\\
(2) All $D_i$ intersect along the boundary ${\bar Y}_{g-1}-Y_{g-1}.$}
\vspace{0.1in}\\
We have several {\it natural} questions.
\vspace{0.1in}\\
\noindent{\bf Problem 6.2.}\quad Describe ${\bar Y}_{g-1}$ and ${\bar Y}_{g-1}-Y_{g-1}$ explicitly in terms of
Jacobi forms. More generally,{\it describe} ${\bar Y}_r$ and ${\bar Y}_r-Y_r$
when $Y_r:=\G_r(k)\ltimes H_{\Bbb Z}^{(r,k)}\ba
H_r\times {\Bbb C}^{(r,k)}\ (\,1\leq r\leq g\,). $
\vspace{0.1in}\\
\noindent{\bf Problem 6.3.}\quad Describe the field of meromorphic functions on ${\bar Y}_{g-1}$ or
${\bar Y}_r.$
\vspace{0.1in}\\
\noindent{\bf Problem 6.4.}\quad Can any $\G_{g-1}^J(k)$-invariant or $\G_r^J(k)$-invariant
meromorphic function on $Y_{g-1}$ or
$Y_r$ be expressed by a quotient of two Jacobi forms of the same weight and
index?
\vspace{0.1in}\\
\noindent{\bf 7. The Boundary of the Satake Compactification}
\setcounter{equation}{0}
\renewcommand{\theequation}{7.\arabic{equation}}
\vspace{0.1in}\\
\indent Let $\G$ be a discrete subgroup of $Sp(g,\BQ)$ which is commensurable
with $\G_g.$ We denote by $M_k(\G)$ the complex vector space consisting of
Siegel modular forms of weight $k$ with respect to $\G\, (\,k\in \BZ\,).$
These vector spaces generate a positively graded ring
$$M(\G):\,=\,\oplus_{k\geq 0}\,M_k(\G)$$
which are integrally closed and of finite type over $M_0(\G)=\BC.$ The
projective variety $A_{g,\G}^*$ associated with $M(\G)$ contains a Zariski
open subset which is complex analytically isomorphic to $A_{g,\G}:=
\G\ba H_g.$ In addition, the boundary $\partial A_{g,\G}^*:=\,A_{g,\G}^*
-A_{g,\G}$ is a disjoint union of a finite number of rational boundary
components of $H_g.$
\vspace{0.1in}\\
\indent From now on, we let $\G:=\,\G_g(k)$ be the principal congruence subgroup
of $\G_g$ of level $k$. We write $g=p+q$ for $0\leq p < g.$ We write an
element $Z$ of $H_g$ as
$$\begin{pmatrix} \tau & W \\
^t\!W & T \end{pmatrix},\ \ \ \tau\in H_p,\ W\in \BC^{(p,q)},\ \ T\in H_q,$$
or simply $Z=(\tau,W,T).$ The Siegel operator $\Phi:M(\G_g(k))\lrt
M(\G_p(k))$ defined by
\begin{equation}
\left(\Phi f\right)(\tau):=\,\lim_{\text{Im}\,T\rightarrow 0}
f\left(\begin{pmatrix} \tau & W\\ * & T\end{pmatrix}\right)\,=\,
\lim_{c\rightarrow 0}
f\left( \begin{pmatrix} \tau & 0 \\ 0 & icE_q \end{pmatrix}\right)
\end{equation}
is a weight-preserving homomorphism which is almost surjective in the sense
that it is surjective for all large weights. Thus we have a canonical
holomorphic embedding $\Phi^*:A_{p,\G_p(k)}^*\lrt A_{g,\G_g(k)}^*.$
We can see that the image of $A_{p,\G_p(k)}=\G_p(k)\ba H_p$ is a
quasi-projective subvariety of $A_{g,\G_g(k)}^*$ and that $Sp(g,\,\BZ/k\BZ)$
acts on $A_{g,\G_g(k)}^*$ as automorphisms. $Sp(g,\,\BZ/k\BZ)$ transforms
$\Phi^*(A_{p,\G_p(k)})$ to its conjugates. Thus we have
$$
\begin{array}{ll}
\partial A_{g,\G_g(k)}^*:&=A_{g,\G_g(k)}^*-A_{g,\G_g(k)}\\
&=\,
\bigcup_{\gamma\in Sp(g,\,\BZ/k\BZ)}\coprod_{l=0}^{g-1}\,\gamma\cdot
\Phi^*(A_{l,\G_l(k)})
\end{array}$$
So in order to investigate the boundary $\partial A_{g,\G_g(k)}^*,$ it is
enough to investigate the boundary points in the image
$\Phi^*(A_{p,\G_p(k)})$ of $A_{p,\G_p(k)}=\G_p(k)\ba H_p$ under $\Phi^*$
for $0\leq p < g.$
\vspace{0.1in}\\
\indent Omitting the detail, we state the following results.\vspace{0.1in}\\
\noindent{\bf Theorem 7.1(Igusa).}\quad{\it Let $\tau_0$ be an element of $H_p.$
Then the analytic local ring ${\mathcal{O}}$ of $A_{g,\G_g(k)}^*$ at the image
point of $\tau_0$ under $\Phi^*$ consists of convergent series of the
following form
$$f(\tau,W,T)=\sum_{\M}\left(\sum_{u}\phi_{\M}(\tau,W\,^tu)\,
e^{\frac{2\pi i\s(\M [u]T)}{k}}\right),\ \ \ \phi_{\M}\in
J_{0,\M}(\G_g(k)),$$
where $\M$ runs over the equivalent classes of
inequivalent half-integral semi-positive symmetric
matrices of degree $q,\ \phi_{\M}$ is a holomorphic function defined on
$V\times \BC^{(q,p)}$ for some open neighborhood $V$ of $\tau_0$ in
$H_p$ and $u$ runs over distinct $\M[u]$ for
$u\in GL(q,\BZ)(k).$}
\vspace{0.1in}\\
\noindent{\bf Theorem 7.2(Igusa).}\quad{\it The ideal $I$ in ${\mathcal{O}}$ associated with
the boundary $\partial A_{g,\G_g(k)}^*\,=\,A_{g,\G_g(k)}^*-A_{g,\G_g(k)}$
consists of convergent series
$$\sum_{\M}\left(\sum_{u}\phi_{\M}(\tau,W\,^tu)\,e^{\frac{2\pi i\s(\M [u]T)}{k}}
\right),\ \ \ \phi_{\M}\in J_{0,\M}(\G_g(k)),$$
where $\M$ runs over inequivalent symmetric positive definite half-integral
matrices of degree $q,\ \phi_{\M}$ is a holomorphic function defined on
$V\times \BC^{(q,p)}$ for some open neighborhood $V$ of $\tau_0$
in $H_p$ and $u$ runs over distinct $\M[u]$ for $u\in
GL(q,\BZ)(k).$}
\vspace{0.1in}\\
\noindent{\bf 8. Singular Jacobi Forms}
\setcounter{equation}{0}
\renewcommand{\theequation}{8.\arabic{equation}}
\vspace{0.1in}\\
\indent In this section, we discuss the notion of singular Jacobi forms. Without
loss of generality we may assume that $\M$ is positive definite. For
simplicity, we consider the case that  $\G$ is the Siegel modular
group $\G_g$ of degree $g.$
\vspace{0.1in}\\
\indent Let $g$ and $h$ be two positive integers. We recall that $\M$ is a symmetric
positive definite, half-integral matrix of degree $h$.
We let
$${\mathcal{P}}_g:=\{ Y\in \BR^{(g,g)}\,\vert\, Y=\,{^tY} > 0\,\} $$
be the open convex cone of positive definite matrices of degree $g$
in the Euclidean space $\BR^{{g(g+1)}\over 2}.$
We define the differential
operator $M_{g,h,\M}$ on ${\mathcal{P}}_g\times \BR^{(h,g)}$ defined by
$$M_{g,h,{\mathcal{M}}}:=\text{det}
\,(Y)\cdot \text{det}\left( {{\partial}\over {\partial Y}}
+{1\over {8\pi}} ^t\!\left( {{\partial}\over {\partial V}}\right)
{\mathcal{M}}^{-1}
\left( {{\partial}\over
{\partial V}}\right) \right), $$
where
$$Y=(y_{\mu\nu})\in {\mathcal{P}}_g,\ \ V=(v_{kl})\in \BR^{(h,g)},\ \
{{\partial}\over {\partial Y}}=\left( {{1+\delta_{\mu\nu}}\over
2}{{\partial}\over {\partial y_{\mu\nu}}}\right)$$ and
$${{\partial}\over {\partial V}}=\left( { {\partial}\over {\partial v_{kl}} }
\right).$$
Yang [85] characterized singular Jacobi forms as follows\,:
\vspace{0.1in}\\
\noindent{\bf Theorem 8.1.}\quad{\it Let $f\in J_{\rho,\M}(\G_g)$ be a Jacobi form of index
$\M$ with respect to a finite dimensional rational representation $\rho$ of
$GL(g,\BC).$ Then the following conditions are equivalent :\\
\indent (1) $f$ is a {\it singular} Jacobi form.\\
\indent (2) $f$ satisfies the differential equation $M_{g,h,\M}f=0.$}
\vspace{0.1in}\\
\noindent{\bf Theorem 8.2.}\quad{\it Let $\rho$ be an irreducible finite dimensional
representation of $GL(g,\BC).$ Then there exists a nonvanishing
{\it singular} Jacobi form in $J_{\rho,\M}(\G_g)$ if and only if
$2k(\rho)< g+h.$ Here $k(\rho)$ denotes the weight of $\rho.$}
\vspace{0.1in}\\
\indent For the proofs of the above theorems we refer to [85], Theorem
4.1 and Theorem 4.5.
\vspace{0.1in}\\
\noindent {\bf Exercise 8.3.}\quad Compute the eigenfunctions and the eigenvalues of
$M_{g,h,\M}$(cf. [85], pp. 2048-2049).
\vspace{0.1in}\\
\indent Now we consider the following group $GL(g,\BR)\ltimes
H_{\BR}^{(g,h)}$ equipped with the multiplication law $$\begin{array}{ll} &\
(A,(\l,\mu,\k))\ast (B,(\l',\mu',\k')) \\ &=(AB,\,(\l
B+\l',\mu\,^t\!B^{-1}+\mu',\k+\k'+\l B\,^t\!\mu'-\mu\,^t\!B^{-1}
\,^t\!\l')),\end{array}$$ where $A,B\in GL(g,\BR)$ and
$(\l,\mu,\k),(\l',\mu',\k')\in H_{\BR}^{(g,h)}.$ We observe that
$GL(g,\BR)$ acts on $H_{\BR}^{(g,h)}$ on the right as
automorphisms. And we have the canonical action of
$GL(g,\BR)\ltimes H_{\BR}^{(g,h)}$ on ${\mathcal{P}}_g\times
\BR^{(h,g)}$ defined by
\begin{equation}
(A,(\l,\mu,\k))\circ (Y,V):=(AY\,^t\!A,\,(V+\l Y+\mu)\,^tA),
\end{equation}
where $A\in GL(g,\BR), (\l,\mu,\k)\in H_{\BR}^{(g,h)}$ and $(Y,V)\in {\mathcal{P}}_g\times
\BR^{(h,g)}.$
\vspace{0.1in}\\
\noindent{\bf Lemma 8.4.}\quad {\it $M_{g,h,\M}$ is invariant under the action of
$GL(g,\BR)\ltimes \left\{\,(0,\mu,0)\,\vert\\mu\in{\BR}^{(h,g)}\,\right\}.$}
\vspace{0.05in}\\
\noindent{\it Proof.}\quad It follows immediately from the direct calculation.
\vspace{0.1in}\\
\indent We have the following natural questions.
\vspace{0.1in}\\
\noindent{\bf Problem 8.5.}\quad Develope the invariant theory for the action of
$GL(g,\BR)\ltimes H_{\BR}^{(g,h)}$ on ${\mathcal{P}}_g\times \BR^{(h,g)}.$
\vspace{0.1in}\\
\noindent{\bf Problem 8.6.}\quad Discuss the application of the theory of singular
Jacobi forms to the geometry of the universal abelian variety as that of
singular modular forms to the geometry of the Siegel modular
variety\,(\,see Appendix B\,).
\vspace{0.1in}\\
\noindent{\bf 9. The Siegel-Jacobi Operator}
\setcounter{equation}{0}
\renewcommand{\theequation}{9.\arabic{equation}}
\vspace{0.1in}\\
\indent  Let $\rho$ and $\M$ be the same as in the previous sections. For positive
integers $r$ and $g$ with $r<g,$ we let $\rho^{(r)}:GL(r,\BC)\lrt
GL(V_{\rho})$ be a rational representation of $GL(r,\BC)$ defined by
$$\rho^{(r)}(a)v:=\rho\left( \begin{pmatrix} a & 0\\ 0 & E_{g-r}\end{pmatrix}
\right)v,\ \ \ a\in GL(r,\BC),\ \ v\in V_{\rho}.$$
The Siegel-Jacobi operator $\Psi_{g,r}:J_{\rho,\M}(\G_g)\lrt
J_{\rho^{(r)},\M}(\G_r)$ is defined by
\begin{equation}
(\Psi_{g,r}f)(Z,W):=\lim_{t\rightarrow \infty}\,f\left(
\begin{pmatrix} Z & 0\\ 0 & itE_{g-r}\end{pmatrix},(W,0)\right),
\end{equation}
where $f\in J_{\rho,\M}(\G_g),\ Z\in H_r$ and $W\in \BC^{(h,r)}.$
It is easy to check that the above limit always exists and the
Siegel-Jacobi operator is a linear mapping. Let $V_{\rho}^{(r)}$ be the
subspace of $V_{\rho}$ spanned by the values
$\{\,(\Psi_{g,r}f)(Z,W)\,\vert\ f\in J_{\rho,\M}(\G_g),\ (Z,W)\in
H_r\times \BC^{(h,r)}\,\}.$ Then $V_{\rho}^{(r)}$ is invariant under the
action of the group
$$\left\{\,\begin{pmatrix} a & 0\\ 0 & E_{g-r}\end{pmatrix}\,:\ a\in GL(r,\BC)\
\right\}\cong GL(r,\BC).$$
We can show that if $V_{\rho}^{(r)}\neq 0$ and $(\rho,V_{\rho})$ is
irreducible, then $(\rho^{(r)},V_{\rho}^{(r)})$ is also irreducible.
\vspace{0.1in}\\
\noindent{\bf Theorem 9.1.}\quad{\it The action of the Siegel-Jacobi operator is compatible
with that of that of the Hecke operator.}
\vspace{0.1in}\\
\indent We refer to [83] for a precise detail on the Hecke operators and
the proof of the above theorem.
\vspace{0.1in}\\
\noindent{\bf Problem 9.2.}\quad Discuss the injectivity, surjectivity and bijectivity
of the Siegel-Jacobi operator.
\vspace{0.1in}\\
\indent This problem was partially discussed by Yang [83] and Kramer [40] in the
special cases. For instance, Kramer [40] showed that if $g$ is arbitrary, $h=1$ and $\rho:GL(g,\BC)
\lrt \BC^{\times}$ is a one-dimensional representation of $GL(g,\BC)$
defined by $\rho(a):=(\text{det}\,(a))^k$ for some $k\in \BZ^+,$ then
the Siegel-Jacobi operator $$\Psi_{g,g-1}:J_{k,m}(\G_g)\lrt J_{k,m}(\G_{g-1})$$
is surjective for $k\gg m\gg 0.$
\vspace{0.1in}\\
\noindent{\bf Theorem 9.3.}\quad{\it Let $1\leq r\leq g-1$ and let $\rho$ be an irreducible
finite dimensional representation of $GL(g,\BC).$ Assume that
$k(\rho)> g+r+\text{rank}\,(\M)+1$ and that $k$ is even. Then
$$J_{\rho^{(r)},\M}^{\text{cusp}}(\G_r)\subset
\Psi_{g,r}(J_{\rho,\M}(\G_g)).$$
Here $J_{\rho^{(r)},\M}^{\text{cusp}}(\G_r)$ denotes the subspace consisting
of all cuspidal Jacobi forms in $J_{\rho^{(r)},\M}(\G_r).$}
\vspace{0.05in}\\
\noindent{\it Idea of Proof.}\quad For each $f\in J_{\rho^{(r)},\M}^{\text{cusp}}
(\G_r),$ we can show by a direct computation that
$$\Psi_{g,r}(E_{\rho,\M}^{(g)}(Z,W;f))=f,$$
where $E_{\rho,\M}^{(g)}(Z,W;f)$ is the Eisenstein series of Klingen's
type associated with a cusp form $f.$ For a precise detail, we refer to [94].
\vspace{0.1in}\\
\noindent{\bf Remark 9.4.}\quad Dulinski [17] decomposed the vector space
$J_{k,\M}(\G_g)\,(k\in \BZ^+)$ into a direct sum of certain subspaces
by calculating the action of the Siegel-Jacobi operator on
Eisenstein series of Klingen's type explicitly.
\vspace{0.1in}\\
\indent For two positive integers $r$ and $g$ with $r\leq g-1,$ we consider the
bigraded ring
$$J_{\ast,\ast}^{(r)}(\ell):=\oplus_{k=0}^{\infty}\oplus_{\M}
J_{k,\M}(\G_r(\ell))$$
and
$$M_{\ast}^{(r)}(\ell):=\oplus_{k=0}^{\infty}\,J_{k,0}(\G_r(\ell))
=\oplus_{k=0}^{\infty}[\G_r(\ell),k],$$
where $\G_r(\ell)$ denotes the principal congruence subgroup of $\G_r$
of level $\ell$ and $\M$ runs over the set of all symmetric semi-positive
half-integral matrices of degree $h$.
Let
$$\Psi_{r,r-1,\ell}:J_{k,\M}(\G_r(\ell))\lrt J_{k,\M}(\G_{r-1}(\ell))$$
be the Siegel-Jacobi operator defined by (9.1).
\vspace{0.1in}\\
\noindent{\bf Problem 9.5.}\quad Investigate $\text{Proj}\,J_{\ast,\ast}^{(r)}(\ell)$ over
$M_{\ast}^{(r)}(\ell)$ and the quotient space
$$Y_r(\ell):=
({{\G_r(\ell)\ltimes (\ell\BZ)^2})\ba ({H_r\ltimes \BC^r}})$$ for
$1\leq r\leq g-1.$\vspace{0.1in}\\
The difficulty to this problem comes from the following facts (A) and
(B)\,:\\
\indent (A) $J_{\ast,\ast}^{(r)}(\ell)$ is not finitely generated over
$M_{\ast}^{(r)}(\ell).$\\
\indent (B) $J_{k,\M}^{\text{cusp}}(\G_r(\ell))\neq \text{ker}\,\Psi_{r,r-1,\ell}$
in general.
\vspace{0.1in}\\
\indent These are the facts different from the theory of Siegel modular forms.
We remark that Runge([54], pp. 190-194) discussed some parts about the above problem.
\vspace{0.1in}\\
\noindent{\bf 10. Invariant Metrics on the Siegel-Jacobi Space}
\setcounter{equation}{0}
\renewcommand{\theequation}{10.\arabic{equation}}
\vspace{0.1in}\\
\def\HC{H_1\times {\mathbb C}}
\def\BD{\mathbb D}
\def\BH{\mathbb H}
\def\BR{\mathbb R}
\def\BC{\mathbb C}
\def\lrt{\longrightarrow}
\def\lmt{\longmapsto}
\def\CX{{\mathcal X}}
\def\td{\bigtriangledown}
\def\pdx{ {{\partial}\over{\partial x}} }
\def\pdy{ {{\partial}\over{\partial y}} }
\def\pdu{ {{\partial}\over{\partial u}} }
\def\pdv{ {{\partial}\over{\partial v}} }
\def\PZ{ {{\partial}\over {\partial Z}} }
\def\PW{ {{\partial}\over {\partial W}} }
\def\PZB{ {{\partial}\over {\partial{\overline Z}}} }
\def\PWB{ {{\partial}\over {\partial{\overline W}}} }
\def\PX{ {{\partial}\over{\partial X}} }
\def\PY{ {{\partial}\over {\partial Y}} }
\def\PU{ {{\partial}\over{\partial U}} }
\def\PV{ {{\partial}\over{\partial V}} }
\def\th{\theta}
\def\l{\lambda}
\def\k{\kappa}
\def\G{\Gamma}
\def\s{\sigma}
\def\g{\gamma}
\def\ddx{{{\partial^2}\over{\partial x^2}}}
\def\ddy{{{\partial^2}\over{\partial y^2}}}
\def\ddu{{{\partial^2}\over{\partial u^2}}}
\def\ddv{{{\partial^2}\over{\partial v^2}}}
\def\px{{{\partial}\over{\partial x}}}
\def\py{{{\partial}\over{\partial y}}}
\def\pu{{{\partial}\over{\partial u}}}
\def\pv{{{\partial}\over{\partial v}}}
\def\pxu{{{\partial^2}\over{\partial x\partial u}}}
\def\pyv{{{\partial^2}\over{\partial y\partial v}}}
\def\DSPR{{\mathbb D}(\SPR)}
\def\dx{{{\partial}\over{\partial x}}}
\def\dy{{{\partial}\over{\partial y}}}
\def\du{{{\partial}\over{\partial u}}}
\def\dv{{{\partial}\over{\partial v}}}
\def\bz{d{\overline Z}}
\def\bw{d{\overline W}}
\def\SJ{H_g\times {\mathbb C}^{(h,g)}}
\def\tr{\triangledown}
\def\Hnm{ H_{g,h}}
\def\Hn{H_n}
\def\Cmn{{\mathbb C}^{(h,g)}}
\indent For a brevity, we write $\Hnm:=\SJ.$ For a coordinate
$(Z,W)\in\Hnm$ with $Z=(z_{\mu\nu})\in {\Bbb H}_g$ and
$W=(w_{kl})\in \Cmn,$ we put $$\begin{array}{ll} Z\,=&\,X\,+\,iY,\quad\ \
X\,=\,(x_{\mu\nu}),\quad\ \ Y\,=\,(y_{\mu\nu}) \ \ \text{real},\\
W\,=&U\,+\,iV,\quad\ \ U\,=\,(u_{kl}),\quad\ \ V\,=\,(v_{kl})\ \
\text{real},\\ dZ\,=&\,(dz_{\mu\nu}),\quad\ \
dX\,=\,(dx_{\mu\nu}),\quad\ \ dY\,=\,(dy_{\mu\nu}),\\
dW\,=&\,(dw_{kl}),\quad\ \ dU\,=\,(du_{kl}),\quad\ \
dV\,=\,(dv_{kl}),
\end{array}$$
$$\begin{array}{ll} \PZ\,=\,&\left(\, { {1+\delta_{\mu\nu}} \over 2}\, {
{\partial}\over {\partial z_{\mu\nu}} } \,\right),\quad
\PZB\,=\,\left(\, { {1+\delta_{\mu\nu}}\over 2} \, {
{\partial}\over {\partial {\overline z}_{\mu\nu} }  } \,\right),\\
\PX\,=\,&\left(\, { {1+\delta_{\mu\nu}}\over 2}\, {
{\partial}\over {\partial x_{\mu\nu} } } \,\right),\quad
\PY\,=\,\left(\, { {1+\delta_{\mu\nu}}\over 2}\, { {\partial}\over
{\partial y_{\mu\nu} }  } \,\right),
\end{array}$$
$$\PW:=\begin{pmatrix} { {\partial}\over{\partial w_{11}} } & \hdots & {
{\partial}\over{\partial w_{h1}} }\\ \vdots&\ddots&\vdots\\ {
{\partial}\over{\partial w_{1g}} }&\hdots &{ {\partial}\over
{\partial w_{hg}} } \end{pmatrix},\quad \PWB:=\begin{pmatrix} {
{\partial}\over{\partial {\overline w}_{11} }   }& \hdots&{
{\partial}\over{\partial {\overline w}_{h1} }  }\\
\vdots&\ddots&\vdots\\ { {\partial}\over{\partial{\overline
w}_{1g} }  }&\hdots & { {\partial}\over{\partial{\overline w}_{hg}
}  }
\end{pmatrix},$$
$$\PU:= \begin{pmatrix} { {\partial}\over{\partial u_{11}} }&\hdots& {
{\partial}\over{\partial u_{h1}} } \\ \vdots &\ddots &\vdots\\ {
{\partial}\over{\partial u_{1g}} } &\hdots & {
{\partial}\over{\partial u_{hg}} }\end{pmatrix},\quad \PV:= \begin{pmatrix}
{ {\partial}\over{\partial v_{11}} } &\hdots & {
{\partial}\over{\partial v_{h1}} }\\ \vdots &\ddots &\vdots\\ {
{\partial}\over{\partial v_{1g}} } &\hdots & {
{\partial}\over{\partial v_{hg}} }\end{pmatrix}.$$ We let
$$T_g:=\,\left\{\,z\in \BC^{(g,g)}\,|\ z=\,^tz\ \right\}$$ be the
vector space of all $g\times g$ complex {\it symmetric} matrices.
The unitary group $K:=U(g)$ of degree $g$ acts on the complex
vector space $T_g\times\Cmn$ by
\begin{equation}
k\cdot(z,w):=\,(\,k\,z\,^tk,\,w\,^tk\,),\quad\ \ k\in U(g),\ z\in T_g,\ w\in \Cmn.
\end{equation}
Then this action induces naturally
the action $\rho$ of $U(g)$ on the polynomial algebra
$\text{Pol}_{h,g}:=\,\text{Pol}\,(T_g\times\Cmn)$. We denote by
$\text{Pol}_{h,g}^K$ the subalgebra of $\text{Pol}_{h,g}$
consisting of all $K$-invariants of the action $\rho$ of
$K:=U(g).$ We also denote by $\BD(\Hnm)$ the algebra of all
differential operators on $\Hnm$ which is invariant under the
action (2.1) of the Jacobi group $G_{g,h}^J$. Then we can show
that there exists a natural linear bijection
\begin{equation}
\Phi:\,\text{Pol}^K_{h,g}\lrt \BD(\Hnm)
\end{equation}
of $\text{Pol}^K_{h,g}$ onto $\BD(\Hnm).$
\vspace{0.1in}\\
\noindent{\bf Theorem 10.1.}\quad{\it The algebra $\BD(\Hnm)$ is generated by the
images under the mapping $\Phi$ of the following invariants}
\vspace{0.1in}\\
$\quad\text{(I1)}\ \quad\quad p_j(z,w):=\,\s((z{\bar
z})^j),\quad\quad 1\leq j\leq g,$
\vspace{0.1in}\\
$\quad\text{(I2)}\quad\quad \psi_k^{(1)}(z,w):=\,(w\,^t{\bar
w})_{kk},\quad\quad 1\leq k\leq h,$
\vspace{0.1in}\\
$\quad\text{(I3)}\quad\quad \psi_{kp}^{(2)}(z,w):=$ Re$(w\,^t{\bar w})_{kp}, \quad\quad 1\leq k<p\leq h,$
\vspace{0.1in}\\
$\quad\text{(I4)}\quad\quad \psi_{kp}^{(2)}(z,w):=$ Im$(w\,^t{\bar w})_{kp}, \quad\quad 1\leq k<p\leq h,$
\vspace{0.1in}\\
$\quad\text{(I5)}\quad\quad f_{kp}^{(1)}(z,w):=$ Re$(w{\bar z}\,^tw)_{kp}, \quad\quad 1\leq k\leq p\leq
h$
\vspace{0.1in}\\
\noindent{\it and}
\vspace{0.1in}\\
$\quad\text{(I6)}\quad\quad f_{kp}^{(2)}(z,w):=$ Im$(w{\bar z}\,^tw)_{kp}, \quad\quad 1\leq k\leq p\leq
h.$\\

In particular, $\BD(\Hnm)$ is not commutative.
\vspace{0.1in}\\
\noindent{\bf Theorem\ 10.1'.}\quad {\it The algebra $\BD( H_{1,1})$ is generated by
the following differential operators $$\begin{array}{ll}
D:=&\,y^2\,\left(\,{{\partial^2}\over {\partial x^2}}\,+\,
{{\partial^2}\over {\partial y^2}}\,\right)
\,+\,v^2\,\left(\,\ddu\,+\,\ddv\,\right) \\ &\ \
+\,2\,y\,v\,\left(\,\pxu\,+\,\pyv\,\right),
\end{array}$$
$$\Psi:=\,y\,\left(\,{{\partial^2}\over {\partial u^2}}\,+\,
{{\partial^2}\over {\partial v^2}}\,\right),$$
$$D_1:=\,2\,y^2\,{{\partial^3}\over {\partial x\partial u
\partial v}}\,-\,y^2\,{{\partial}\over{\partial y}}
\left(\,{{\partial^2}\over{\partial u^2}}\,-\,
{{\partial^2}\over{\partial v^2}}\,\right)\,+\,\left(\,
v\,{{\partial}\over{\partial v}}\,+\,1\,\right)\Psi $$ and
$$D_2:=\,y^2\,{{\partial}\over{\partial x}}\left(\,
{{\partial^2}\over{\partial v^2}}\,-\,{{\partial^2}\over {\partial
u^2}}\,\right)\,-\,2\,y^2\,{{\partial^3}\over{\partial y\partial u
\partial v}}\,-\,v\,{{\partial}\over{\partial u}}\Psi,
$$ where $\tau=x+iy$ and $z=u+iv$ with real variables $x,y,u,v.$
Moreover, we have $$\begin{array}{ll} [D,\,\Psi]:=D\Psi-\Psi D\,=&\,
2\,y^2\,\dy\left(\,\ddu\,-\,\ddv\,\right)\,-\,
4\,y^2\,{{\partial^3}\over{\partial x\partial u\partial v}}\\ &\ \
\,-\,2\,\left(\,v\,\dv\Psi\,+\,\Psi\,\right).
\end{array}  $$
In particular, the algebra $\BD(H_{1,1})$ is not commutative.}
\vspace{0.1in}\\
\noindent {\bf Theorem 10.2.}\quad{\it The following metric
\begin{equation}
\begin{array}{ll}
ds_{g,h}^2:=& \s\left(Y^{-1}dZ\,Y^{-1}d{\overline Z}\right)\,+ \,\s\left(Y^{-1}\,^tV\,V\,Y^{-1}dZ\,Y^{-1} \bz\right)\\
            & \,\,\,\,+\,\s\left(Y^{-1}\,^t(dW)\,\bw\right)\\
            & \s\left(Y^{-1}dZ\,Y^{-1}\,^t(\bw)\,V\, +\,Y^{-1}\bz\,Y^{-1}\,^t(dW)\,V\right)
\end{array}
\end{equation}
is a Riemannian metric on the Siegel-Jacobi space $\Hnm$ which is
invariant under the action (1.2) of the Jacobi group $G^J_{g,h}$.
Also the above metric is a K{\" a}hler metric. The
Laplace-Beltrami operator $\Delta_{g,h}$ of the Siegel-Jacobi
space ($\Hnm,\,ds_{g,h}^2)$ is given by
\begin{equation}
\begin{array}{ll}
\Delta_{h,g}\,=\,& 4\,\s\left(\,Y\,\PZ\,Y\,\PZB\,\right)\,+\,
4\,\s\left(\, Y\,\PW\,^t\left( \PWB\right)\,\right)\\ &\ \ \ \
+\,4\,\s\left(\,\PW\,V\,\PWB\,V\,\right)\\
&\ \
+\,4\,\s\left(\,\PZ\,Y\,\PWB\,V\,+\,\PZB\,Y\,\PW\,V\,\right).
\end{array}
\end{equation}
\def\w{\wedge}
The following differential form
$$dv:=\,\left(\,\text{det}\,Y\,\right)^{-(g+h+1)}[dX]\w [dY]\w
[dU]\w [dV]$$ is a $G_{g,h}^J$-invariant volume element on $\Hnm$,
where $$[dX]:=\,\w_{\mu\leq\nu}dx_{\mu\nu},\quad
[dY]:=\,\w_{\mu\leq\nu} dy_{\mu\nu},\quad
[dU]:=\,\w_{k,l}du_{kl}\quad\text{and}\quad
[dV]:=\,\w_{k,l}dv_{k,l}.$$}
\noindent{\bf Theorem 10.3.}\quad{\it The automorphism group of $\Hnm$ is
isomorphic to the group $Sp(g,\BR)\ltimes
\left(\,\BR^{(h,g)}\times\BR^{(h,g)}\,\right)$ equipped with the
multiplication $$(M,(\l,\mu))\cdot
(M',(\l',\mu)):=\,(MM',\,({\tilde{\l}}+\l', {\tilde{\mu}}+
\mu')),$$ where $M,M'\in Sp(g,\BR),\ \l,\mu\in \BR^{(h,g)}$ and
$({\tilde{\l}},{\tilde{\mu}}):=(\l,\mu)M'.$}
\vspace{0.1in}\\
\noindent{\bf Theorem 10.4.}\quad{\it The scalar curvature of the Siegel-Jacobi
space $(\HC,\,ds^2)$ is $-3$}.
\vspace{0.1in}\\
\indent We note that according to Theorem 2, the metric $ds^2$ is given by
\begin{equation}
\begin{array}{ll}
ds^2:=ds_{1,1}^2=&{{y+v^2}\over {y^3}}\,(dx^2+dy^2)\,+\,{\frac 1y}\,(du^2+dv^2)\\
                 &\ \ \ \ \ \ -{{2v}\over{y^2}}\,(dxdu+dydv)
\end{array}
\end{equation}
on $\HC$ which is invariant under the action (2.1) of the Jacobi
group $G_{1,1}^J=SL(2,\BR)\ltimes H_{\BR}^{(1,1)}$, where
$z=x+iy\in H_1$ and $w=u+iv\in \BC$ with $x,y,u,v$ real
coordinates.
\vspace{0.1in}\\
\noindent{\bf Remark 10.5.}\quad The Poincar{\'e} upper half plane $H_1$ is a
two dimensional Riemannian manifold with the Poincar{\'e} metric
$$ds_0^2:=\,{{dx^2\,+\,dy^2}\over {y^2}},\ \ \ z=x+iy\in H_1 \
\text{with}\ x,y\ \text{real}.$$ It is easy to see that the
Gaussian curvature is $-1$ everywhere and $H_1$ is an {\it
Einstein\ manifold}. In fact, if we denote by $S_0(X,Y)$ the Ricci
curvature of $(H_1,\,ds_0^2),$ then we have
$$S_0(X,Y)\,=\,-g_0(X,Y)\ \ \ \ \text{for\ all}\ X,Y\in
\CX(H_1),$$ where ${\mathcal {X}}(H_1)$ denotes the algebra of all
smooth vector fields on $H_1$ and $g_0(X,Y)$ is the inner product
on the tangent bundle $T(H_1)$ induced by the Poincar{\'e} metric
$ds_0^2.$ But the Siegel-Jacobi space $\HC$ is {\it not} an
Einstein manifold. Indeed, if we denote by $S(X,Y)$ the Ricci
curvature of $(\HC,\,ds^2)$ and $E_1:=\,{ {\partial}\over
{\partial x} }$, we can see without difficulty that there does not
exist a constant $c$ such that
$$S(E_1,E_1)\,=\,c\,g(E_1,E_1)\,=\,c\,g_{11}\,=\,c\, {
{y+v^2}\over {y^3} },.$$ where $g=(g_{ij})$ is the inner product
on the tangent bundle $T(H_1\times\BC)$ induced by the metric
(10.5).
\vspace{0.1in}\\
\indent Now we will introduce the notion of {\it Maass}-{\it Jacobi\
forms}.
\vspace{0.1in}\\\noindent
{\bf Definition\ 10.6.}\quad A smooth function $f:\,\Hnm\lrt \BC$ is
called a {\it Maass}-{\it Jacobi\ form} on $\Hnm$ if $f$ satisfies
the following conditions (MJ1)-(MJ3) :
\vspace{0.1in}\\
\indent (MJ1) $f$ is invariant under $\G_{g,h}^J:=\,\G_g\ltimes
H_{\Bbb Z}^{(g,h)}.$
\vspace{0.1in}\\
\indent (MJ2) $f$ is an eigenfunction of the Lapalce-Beltrami
operator $\Delta_{n,m}$.
\vspace{0.1in}\\
\indent (MJ3) $f$ has a polynomial growth.
\vspace{0.1in}\\
\noindent Here $\G_g:=Sp(g,{\Bbb Z})$ denotes the Siegel modular group of
degree $g$ and and $$H_{\Bbb Z}^{(g,h)}:=\,\left\{\,(\l,\mu;\k)\in
H_{\BR}^{(g,h)}\,|\ \l,\mu,\k\ \text{integral}\ \right\}.$$
\vspace{0.1in}\\
\noindent For more details on Maass-Jacobi forms in the case $g=h=1,$ we
refer to [89].
\vspace{0.1in}\\
\noindent{\bf 11. Final Remarks}
\setcounter{equation}{0}
\renewcommand{\theequation}{11.\arabic{equation}}
\vspace{0.1in}\\
In [32] and [34], Gritsenko, Hulek and Sankaran gave applications of
Jacobi forms of degree $1$ in the study of the moduli space of abelian
surfaces with a certain polarization. We refer to
[7],[9],[11],[61],[62]
for the representation theory of the Jacobi group.
\vspace{0.1in}\\
\noindent{\bf Appendix A. Subvarieties of the Siegel Modular Variety}
\setcounter{equation}{0}
\renewcommand{\theequation}{10.\arabic{equation}}
\vspace{0.1in}\\
\indent Here we assume that the ground field is the complex number field $\BC.$
\vspace{0.1in}\\
\noindent{\bf Definition A.1.}\quad A nonsingular variety $X$ is said to be
{\it rational} if $X$ is birational to a projective space $\BP^n(\BC)$ for
some integer $n$. A nonsingular variety $X$ is said
to be {\it stably\ rational} if $X\times \BP^k(\BC)$ is birational to $\BP^N(\BC)$ for certain nonnegative
integers $k$ and $N$. A nonsingular variety $X$ is called {\it unirational} if there
exists a dominant rational map $\varphi:\BP^n(\BC)\lrt X$ for a certain
positive integer $n$, equivalently if the function field $\BC(X)$ of $X$
can be embedded in a purely transcendental extension $\BC(z_1,\cdots,z_n)$ of $\BC.$
\vspace{0.1in}\\
\noindent{\bf Remarks A.2.}\quad (1) It is easy to see that the rationality implies
the stably rationality and that the stably rationality implies the
unirationality.
\vspace{0.05in}\\
\noindent (2) If $X$ is a Riemann surface or a complex surface, then the notions of
rationality, stably rationality and unirationality are equivalent one
another.
\vspace{0.05in}\\
\noindent (3) Griffiths and Clemens(cf. Ann. of Math. 95(1972), 281-356)
showed that most of cubic threefolds in $\BP^4(\BC)$ are unirational but {\it not} rational.
\vspace{0.1in}\\
The following natural questions arise :
\vspace{0.1in}\\
{\sc Question 1.}\quad Is a stably rational variety {\it rational}\,? Indeed,
the question was raised by Bogomolov.
\vspace{0.1in}\\
{\sc Question 2.}\quad Is a general hypersurface $X\subset \BP^{n+1}(\BC)$ of
degree $d\leq n+1$ {\it unirational}\,?
\vspace{0.1in}\\
\noindent{\bf Definition A.3.}\quad Let $X$ be a nonsingular variety of dimension $n$
and let $K_X$ be the canonical divisor of $X$. For each positive integer
$m\in \BZ^+$, we define the $m$-{\it genus} $P_m(X)$ of $X$ by
$$P_m(X):=\text{dim}_{\BC}\,H^0(X,{\mathcal{O}}(mK_X)).$$
The number $p_g(X):=P_1(X)$ is called the {\it geometric\ genus} of $X$.
We let
$$N(X):=\left\{\,m\in \BZ^+\,\vert\, P_m(X)\geq 1\,\right\}.$$
For the present, we assume that $N(X)$ is nonempty. For each $m\in N(X),$
we let $\left\{ \phi_0,\cdots,\phi_{N_m}\right\}$ be a basis of the vector
space $H^0(X,{\mathcal{O}}(mK_X)).$ Then we have the mapping $\Phi_{mK_X}\,:
X\lrt \BP^{N_m}(\BC)$ by
$$\Phi_{mK_X}(z):=(\phi_0(z):\cdots:\phi_{N^m}(z)),\ \ \ z\in X.$$
We define the {\it Kodaira\ dimension} $\kappa(X)$ of $X$ by
$$\kappa(X):=\text{max}\,\left\{\,\text{dim}_{\BC}\,\Phi_{mK_X}(X)\,\vert\,\,
m\in N(X)\,\right\}.$$
If $N(X)$ is empty, we put $\kappa(X):=-\infty.$ Obviously $\kappa(X)\leq
\text{dim}_{\BC}\,X.$ A nonsingular variety $X$
is said to be {\it of\ general\
type} if $\kappa(X)=\text{dim}_{\BC}X.$
A singular variety $Y$ in general is said to be rational, stably rational,
unirational or of general type if any nonsingular model $X$ of $Y$ is
rational, stably rational, unirational or of general type respectively.
We define
$$P_m(Y):=P_m(X)\ \ \ \ \text{and}\ \ \ \ \kappa(Y):=\kappa(X).$$
A variety $Y$ of dimension $n$ is said to be {\it of\ logarithmic\
general\ type} if there exists a smooth compactification ${\tilde Y}$ of
$Y$ such that $D:={\tilde Y}-Y$ is a divisor with normal crossings only and
the transcendence degree of the logarithmic canonical ring
$$\oplus_{m=0}^{\infty}\,H^0({\tilde Y},\,m(K_{\tilde Y}+[D]))$$
is $n+1$, i.e., the {\it logarithmic\ Kodaira\ dimension} of $Y$ is $n$.
We observe that the notion of being of logarithmic general type is weaker
than that of being of general type.
\vspace{0.1in}\\
\indent Let $A_g:=\G_g\ba H_g$ be the Siegel modular variety of degree $g$, that is,
the moduli space of principally polarized abelian varieties of dimension
$g$. So far it has been proved that $A_g$ is of general type for $g\geq 7.$
At first Freitag [24] proved this fact when $g$ is a multiple of $24$.
Tai [60] proved this for $g\geq 9$ and Mumford [46] proved this fact for
$g\geq 7.$ On the other hand, $A_g$ is known to be unirational for
$g\leq 5\,:$ Donagi [16] for $g=5,$ Clemens [15] for $g=4$ and classical for
$g\leq 3.$ For $g=3,$ using the moduli theory of curves, Riemann [53],
Weber [65] and Frobenius [28] showed that $A_3(2):=\G_3(2)\ba H_3$ is a
rational variety and moreover gave $6$ generators of the modular function
field $K(\G_3(2))$ written explicitly in terms of derivatives of odd theta
functions at the origin. So $A_3$ is a unirational variety with a Galois
covering of a rational variety of degree $[\G_3:\G_3(2)]=1,451,520.$ Here
$\G_3(2)$ denotes the principal congruence subgroup of $\G_3$ of level $2.$
Furthermore it was shown that $A_3$ is stably rational(cf. [38], [12]).
For a positive integer $k$, we let $\G_g(k)$ be the principal
congruence subgroup of $\G_g$ of level $k$. Let $A_g(k)$ be the moduli
space of abelian varieties of dimension $g$ with $k$-level structure. It is
classically known that $A_g(k)$ is of logarithmic general type for
$k\geq 3$(cf. [45]). Wang [64] proved that $A_2(k)$ is of general
type for $k\geq 4.$ On the other hand, van der Geer [29] showed that
$A_2(3)$ is rational.
The remaining unsolved problems are summarized as follows\,:
\vspace{0.1in}\\
{\bf Problem 1.}\quad Is $A_3$ rational\,?
\vspace{0.1in}\\
{\bf Problem 2.}\quad Are $A_4,\ A_5$ stably rational or rational\,?
\vspace{0.1in}\\
{\bf Problem 3.}\quad Discuss the (uni)rationality of $A_6.$
\vspace{0.1in}\\
{\bf Problem 4.}\quad What type of varieties are $A_g(k)$ for $g\geq 3$ and
$k\geq 2$\,?
\vspace{0.1in}\\
\indent We already mentioned that $A_g$ is of general type if $g\geq 7.$ It is
natural to ask if the subvarieties of $A_g\,(g\geq 7)$ are of general type,
in particular the subvarieties of $A_g$ of codimension one. Freitag [Fr3]
showed that there exists a certain bound $g_0$ such that for $g\geq g_0,$
each irreducible subvariety of $A_g$ of codimension one is of general type.
Weissauer [Wei2] proved that every irreducible divisor of $A_g$ is of general
type for $g\geq 10.$ Moreover he proved that every subvariety of codimension
$\leq g-13$ in $A_g$ is of general type for $g\geq 13.$ We observe that the
smallest known codimension for which there exist subvarieties of $A_g$ for
large $g$ which are not of general type is $g-1.\ A_1\times A_{g-1}$ is a
subvariety of $A_g$ of codimension $g-1$ which is not of general type.
\vspace{0.1in}\\
\noindent{\bf Remark A.4.}\quad Let $\M_g$ be the coarse moduli space of curves of genus $g$
over $\BC.$ Then $\M_g$ is an analytic subvariety of $A_g$ of dimension
$3g-3.$ It is known that $\M_g$ is unirational for $g\leq 10.$ So the
Kodaira dimension $\kappa(\M_g)$
of $\M_g$ is $-\infty$ for $g\leq 10.$ Harris and Mumford
[H-M] proved that $\M_g$ is of general type  for odd $g$ with $g\geq 25$ and
$\kappa(\M_{23})\geq 0.$
\vspace{0.1in}\\
\noindent{\bf Appendix B. Singular Modular Forms}
\setcounter{equation}{0}
\renewcommand{\theequation}{10.\arabic{equation}}
\vspace{0.1in}\\
\indent Let $\rho$ be a rational representation of $GL(g,\BC)$ on a finite
dimensional complex vector space $V_{\rho}.$ A holomorphic function
$f:H_g\lrt V_{\rho}$ with values in $V_{\rho}$ is called a modular form
of type $\rho$ if it satisfies
$$f(M<Z>)=\rho(CZ+D)f(Z)$$
for all $\begin{pmatrix} A & B\\ C & D \end{pmatrix}\in \G_g$ and $Z\in H_g.$
We denote by $[\G_g,\rho]$ the vector space of all modular forms of type
$\rho.$ A modular form $f\in [\G_g,\rho]$ of type $\rho$ has a Fourier
series $$f(Z)=\sum_{T\geq 0}a(T)e^{2\pi i(TZ)},\ \ \ Z\in H_g,$$
where $T$ runs over the set of all semipositive half-integral symmetric
matrices of degree $g.$ A modular form $f$ of type $\rho$ is said to
be {\it singular} if a Fourier coefficient $a(T)$ vanishes unless
$\text{det}\,(T)=0.$
\vspace{0.1in}\\
\indent Freitag [25] proved that every singular modular form can be written as
a finite linear combination of theta series with harmonic coefficients and
proposed the problem to characterize singular modular forms. Weissauer [66]
gave the following criterion.
\vspace{0.1in}\\
\noindent{\bf Theorem B.1.}\quad{\it Let $\rho$ be an irreducible rational representation of
$GL(g,\BC)$ with its highest weight $(\l_1,\cdots,\l_g).$ Let $f$ be a
modular form of type $\rho.$ Then the following are equivalent\,:
\vspace{0.05in}\\
(a) $f$ is singular.\vspace{0.05in}\\
(b) $2\l_g < g.$}
\vspace{0.1in}\\
\indent Now we describe how the concept of singular modular forms is closely related
to the geometry of the Siegel modular variety. Let $X$ be the Satake
compactification of the Siegel modular variety $A_g=\G_g\ba H_g.$ Then
$A_g$ is embedded in $X$ as a quasiprojective algebraic subvariety of
codimension $g$. Let $X_s$ be the smooth part of $A_g$ and ${\tilde X}$ the
desingularization of $X.$ Without loss of generality, we assume
$X_s\subset {\tilde X}.$ Let $\Omega^p({\tilde X})$\,(resp.\,$\Omega^p(X_s))$
be the space of holomorphic $p$-form on ${\tilde X}$\,(resp.\,$X_s$).
Freitag and Pommerening [27] showed that if $g>1$, then the restriction
map $$\Omega^p({\tilde X})\lrt \Omega^p(X_s)$$
is an isomorphism for $p< \text{dim}_{\BC}\,{\tilde X}={{g(g+1)}\over 2}.$
Since the singular part of $A_g$ is at least codimension $2$ for $g>1,$
we have an isomorphism $$\Omega^p({\tilde X})\cong \Omega^p(H_g)^{\G_g}.$$
Here $\Omega^p(H_g)^{\G_g}$ denotes the space of $\G_g$-invariant holomorphic
$p$-forms on $H_g.$ Let $\text{Sym}^2(\BC^g)$ be the symmetric power of
the canonical representation of $GL(g,\BC)$ on $\BC^n.$ Then we have an
isomorphism
$$\Omega^p(H_g)^{\G_g}\lrt [\G_g,\wedge^p \text{Sym}^2(\BC^g)].$$
\noindent{\bf Theorem B.2([66]).}\quad{\it Let $\rho_{\a}$ be the irreducible representation of
$GL(g,\BC)$ with highest weight $$(g+1,\cdots,g+1,g-\a,\cdots,g-\a)$$
such that $\text{corank}(\rho_{\a})=\a$ for $1\leq \a\leq g.$ If $\a=-1,$
we let $\rho_{\a}=(g+1,\cdots,g+1).$ Then
$$\Omega^p(H_g)^{\G_g}=\begin{cases} [\G_g,\rho_{\a}],\ &\text{if\
$p={{g(g+1)}\over 2}-{{\a(\a+1)}\over 2}$}\\
0, &\text{otherwise}.\end{cases}$$}
\noindent{\bf Remark B.3.}\quad{\it If $2\a>g,$ then any $f\in [\G_g,\rho_{\a}]$ is singular.
Thus if $p<{{g(3g+2)}\over 8},$ then any $\G_g$-invariant holomorphic
$p$-form on $H_g$ can be expressed in terms of vector valued theta series
with harmonic coefficients. It can be shown with a suitable modification
that the just mentioned statement holds for a sufficiently small congruence
subgroup of $\G_g.$}
\vspace{0.1in}\\
\indent Thus the natural question is to ask how to determine the $\G_g$-invariant
holomorphic $p$-forms on $H_g$ for the nonsingular range
$\dfrac{g(3g+2)}{8}\leq p \leq \dfrac{g(g+1)}{2}.$
Weissauer [68] answered the above question for $g=2.$ For $g>2,$ the above question is still open.
It is well know that the vector space of vector valued modular forms of
type $\rho$ is finite dimensional. The computation or the estimate of the
dimension of $\Omega^p(H_g)^{\G_g}$ is interesting because its dimension is
finite even though the quotient space $A_g$ is noncompact.
\vspace{0.1in}\\
\indent Finally we will mention the results due to Weisauer [67].
We let $\G$ be a congruence subgroup of $\G_2.$ According to
Theorem B.2, $\G$-invariant holomorphic forms in $\Omega^2(H_2)^{\G}$ are
corresponded to modular forms of type (3,1). We note that these invariant
holomorphic $2$-forms are contained in the {\it nonsingular\ range}.
And if these modular forms are not cusp forms, they are mapped under the
Siegel $\Phi$-operator to cusp forms of weight $3$ with respect to some
congruence subgroup\,(\,dependent on $\G$\,) of the elliptic modular group.
Since there are finitely many cusps, it is easy to deal with these modular
forms in the adelic version. Observing these facts, he showed that any
$2$-holomorphic form on $\G\ba H_2$ can be expressed in terms of theta
series with harmonic coefficients associated to binary positive definite
quadratic forms. Moreover he showed that $H^2(\G\ba H_2,\BC)$ has a
pure Hodge structure and that the Tate conjecture holds for a suitable
compactification of $\G\ba H_2.$  If $g\geq 3,$ for a congruence subgroup
$\G$ of $\G_g$ it is difficult to compute the cohomology groups
$H^{\ast}(\G\ba H_g,\BC)$ because $\G\ba H_g$ is noncompact and highly
singular. Therefore in order to study their structure, it is natural to
ask if they have pure Hodge structures or mixed Hodge structures.
\vspace{0.2in}\\ 

\footnotesize{
}


\begin{thebibliography}{9}


\bibitem{} A. N. Andrianov ,{\em Modular descent and the Saito-Kurokawa conjecture
}, Invent. Math., {\bf 53}(1979), 267-280.

\bibitem{} T. Arakawa and B. Heim ,{\em Real analytic
Jacobi Eisenstein series and Dirichlet series attached to three
Jacobi forms}, preprint. (MPI 98-66)

\bibitem{} A. Ash, D. Mumford, M. Rapport and Y.-S.
Tai, Smooth compactification of locally symmetric varieties, Math. Sci. Press, Brookline,
(1975).

\bibitem{} R. Berndt, {\em Zur Arithmetik der elliptischen
Funktionenk{\"o}rper h{\"o}herer Stufe}, J. reine angew. Math., {\bf 326}(1981), 79-94.

\bibitem{} R. Bernet, {\em Meromorphic Funktionen auf Mumfords Kompaktifizierung der
universellen elliptischen Kurve $N$-ter Stufe}, J. reine
angew. Math., {\bf 326}(1981), 95-103.

\bibitem{} R. Bernet, {\em Shimuras Reziprozit{\"a}tsgesetz f{\"u}r
den K{\"o}rper der arithmetischen elliptischen Funktionen beliebiger
Stufe }, J. reine angew. Math., {\bf 343}(1983), 123-145.

\bibitem{} R. Bernet, {\em Die Jacobigruppe und die
W{\"a}rmeleitungsgleichung }, Math. Z., {\bf 191}(1986), 351-361.

\bibitem{} R. Bernet, {\em Some remarks on automorphic forms for
the Jacobi group }, IHES/M, (1989).

\bibitem{} R. Bernet, {\em The Continuous Part of $L^2(\G^J\ba G^J)$
for the Jacobi Group}, Abh. Math. Sem. Univ. Hamburg., {\bf 60}(1990),
225-248.

\bibitem{} R. Bernet, {\em On Automorphic Forms for the Jacobi Group
}, Jb. d. Dt. Math.-Verein., {\bf 97}(1995), 1-18.

\bibitem{} R. Berndt and S. B{\"o}cherer, {\em Jacobi Forms
and Discrete Series Representations of the Jacobi Group}, Math.
Z., {\bf 204}(1990), 13-44.

\bibitem{} F. A. Bogomolov and P. I. Katsylo, {\em Rationality
of some quotient varieties}, Math. USSR Sbornik., {\bf 54}(1986),
571-576.

\bibitem{} D. Bump, S. Friedberg and J. Hoffstein, {\em Nonvanishing Theorems for $L$-functions of Modular Forms and
Their Derivatives}, Ann. of Math., {\bf 131}(1990), 53-128.

\bibitem{} J. L. Cardy, {\em Operator Content of Two-dimensional
Conformally Invariant Theories}, Nuclear Physics B, {\bf 270}(1986),
186-204.

\bibitem{} H. Clemens, {\em Double Solids}, Adv. Math., {\bf 47}(1983),
107-230.

\bibitem{} R. Donagi, {\em The unirationality of $A_5$}, Ann. of Math., {\bf 119}(1984),
269-307.

\bibitem{} J. Dulinski, {\em A decomposition theorem for Jacobi
forms}, Math. Ann., {\bf 303}(1995), 473-498.

\bibitem{} M. Eichler, {\em Die Spur der Hecke Operatoren
in gewissen R{\"a}umen von Jacobischen Modulformen}, Abh. Math. Sem. Univ. Hamburg, {\bf 54}(1984),
35-48.

\bibitem{} M. Eichler, {\em Eine neue Klasse von Modulformen
und Modulfunktionen}, Abh. Math. Sem. Univ. Hamburg, {\bf 55}(1985),
53-68.

\bibitem{} M. Eichler, {\em Eine neue Klasse von Modulformen und
Modulfunktionen}, Abh. Math. Sem. Univ. Hamburg, {\bf 57}(1987),
57-68.

\bibitem{} M. Eichler and D. Zagier, The Theory of Jacobi Forms, Progress in Math.,
{\bf 55}, Birkh{\" a}user, Boston-Basel-Stuttgart, (1985).

\bibitem{} G. Faltings and C.-L. Chai, Degeneration of
Abelian Varieties, Springer-Verlag, (1990)

\bibitem{} A. J. Feingold and I. B. Frenkel, {\em A Hyperbolic Kac-Moody Algebra
and the Theory of Siegel Modular Forms of genus $2$}, Math. Ann., {\bf 263}(1983),
87-144.

\bibitem{} E. Freitag, {\em Die Kodairadimension von
K{\"o}rpen automorpher Funktionen}, J. reine angew. Math., {\bf 296}(1977),
162-170.

\bibitem{} E. Freitag, {\em Thetareihnen mit harmonischen
Koeffizienten zur Siegelschen Modulgruppe}, Math. Ann., {\bf 254}(1980),
27-51.

\bibitem{} E. Freitag, {\em Holomorphic tensors on subvarieties of the
Siegel modular variety}, Birkh{\"a}user, Prog. Math. Boston, {\bf 46}(1984),
93-11.

\bibitem{} E. Freitag and K. Pommerening, {\em Regul{\"a}re
Differentialformen des K{\"o}rpers der Siegelschen Modulfunktionen},
J. reine angew. Math., {\bf 331}(1982), 207-220.

\bibitem{} G. Frobenius, {\em {\"U}ber die Beziehungen
zwischen $28$ Doppeltangenten einer eben Curve vierte Ordnung
}, J. reine angew. Math., {\bf 99}(1886), 285-314.

\bibitem{} G. van der Geer, {\em Note on abelian schemes of level
three}, Math. Ann., {\bf 278}(1987), 401-408.

\bibitem{} V. A. Gritsenko, {\em The action of modular operators
on the Fourier-Jacobi coefficients of modular forms}, Math. USSR
Sbornik, {\bf 74}(1984), 237-268.

\bibitem{} B. Gross, W. Kohnen and D. Zagier, {\em Heegner
points and the derivative of $L$-series II}, Math. Ann., {\bf 278}(1987),
497-562.

\bibitem{} V. A. Gritsenko and G. K. Sankaran, {\em Moduli of
abelian surfaces with $(1,p^2)$ polarization}, preprint.

\bibitem{} J. Harris and D. Mumford, {\em On the Kodaira
dimension of the moduli space of curves}, Invent. Math., {\bf 67}(1982),
23-88.

\bibitem{} K. Hulek and G. K. Sankaran, {\em The Kodaira
dimension of certain moduli spaces of abelian surfaces }, Compositio
Math., {\bf 90}(1994), 1-35.

\bibitem{} J. Igusa, {\em A desingularization problem in the
theory of Siegel modular functions}, Math. Ann., {\bf 168}(1967),
228-260.

\bibitem{} C. G. J. Jacobi, {\em Fundamenta nova theoriae
functionum ellipticum}, K{\"o}nigsberg, (1829).

\bibitem{} W. Kohnen, {\em Modular forms of half integral
weight on $\G_0(4)$ }, Math. Ann., {\bf 285}(1980), 249-266.

\bibitem{} J. Kollar and F. O. Schreyer, {\em The moduli of curves
is stably rational for $g\leq 6$ }, Duke Math. J., {\bf 51}(1984),
239-242.

\bibitem{} J. Kramer, {\em A geometrical approach to the
theory of Jacobi forms}, Compositio Math., {\bf 79}(1991), 1-19.

\bibitem{} J. Kramer, {\em An arithmetic theory of Jacobi forms
in higher dimensions }, J. reine angew. Math., {\bf 458}(1995),
157-182.

\bibitem{} N. V. Kuznetsov, {\em A new class of identities for
the Fourier coefficients of modular forms}, Acta Arith., (1975),
505-519.

\bibitem{} H. Maass, {\em {\"U}ber eine Spezialschar von
Modulformen zweiten Grades}, Invent. Math., {\bf 52}(1979),
95-104.

\bibitem{} H. Maass, {\em {\"U}ber eine Spezialschar von
Modulformen zweiten Grades}, Invent. Math., {\bf 53}(1979),
249-253.

\bibitem{} H. Maass, {\em {\"U}ber eine Spezialschar von Modulformen
zweiten Grades}, Invent. Math., {\bf 53}(1979), 255-265.

\bibitem{} D. Mumford, {\em Hirzebruch's proportionality
theorem in the noncompact case}, Invent. Math., {\bf 42}(1977),
239-272.

\bibitem{} D. Mumford, {\em On the Kodaira dimension of the
Siegel modular variety}, Springer Lecture Note, {\bf 997}(1982),
348-375.

\bibitem{} A. Murase, {\em $L$-functions attached to Jacobi
forms of degree $n$. Part I : The Basic Identity}, J. reine angew.
Math., {\bf 401}(1989), 122-156.

\bibitem{} A. Murase, {\em $L$-functions attached to Jacobi forms of
degree $n$. Part II : Functional Equation}, Math. Ann., {\bf 290
}(1991), 247-276.

\bibitem{} A. Murase, {\em On Explicit Formula for the Whittaker-Shintani
Functions on $Sp_2.$}, Abh. Math. Sem. Univ. Hamburg, {\bf 61
}(1991), 153-162.

\bibitem{} A. Murase and T. Sugano, {\em Whittaker-Shintani
Functions on the Symplectic Group of Fourier-Jacobi Type},
Compositio Math., {\bf 79}(1991), 321-349.

\bibitem{} Y. Namikawa, Toroidal compactification of Siegel
spaces, Springer-Verlag, {\bf 812}(1980)

\bibitem{} I. Piateski-Sharpiro, Automorphic Functions
and the Geometry of Classical Domains, Gordan-Breach, New York,
(1966).

\bibitem{} B. Riemann, {\em Zur Theorie der Abel'schen Funktionen
f{\"u}r den $p=3$}, Math. Werke, Teubener, Leipzig, (1876),
456-476.

\bibitem{} B. Runge, {\em Theta functions and Siegel-Jacobi
functions}, Acta Math., {\bf 175}(1995), 165-196.

\bibitem{} I. Satake, {\em Fock Representations and Theta
Functions}, Ann. Math. Studies, {\bf 66}(1971), 393-405.

\bibitem{} I. Satake, Algebraic structures of symmetric domains, Iwanami Shoten
Publ. and Princeton Univ. Press, (1980).

\bibitem{} G. Shimura, {\em On modular forms of half integral
weight }, Ann. of Math., {\bf 97}(1973), 440-481.

\bibitem{} G. Shimura, {\em On certain reciprocity laws for theta
functions and modular forms }, Acta Math., {\bf 141}(1979), 35-71.

\bibitem{} N.-P. Skoruppa, {\em {\"U}ber den Zusammenhang
zwischen Jacobi-Formen und Modulformen halbganzen Gewichts, Dissertation
}, Universit{\"a}t Bonn, (1984).

\bibitem{} Y.-S. Tai, {\em On the Kodaira dimension of the
moduli spaces of abelian varieties}, Invent. Math., {\bf 68}(1982),
425-439.

\bibitem{} K. Takase, {\em A note on automorphic forms}, J. reine angew. Math., {\bf 409}(1990),
138-171.

\bibitem{} K. Takase, {\em On unitary representations of Jacobi groups
}, J. reine angew. Math., {\bf 430}(1992), 130-149.

\bibitem{} W. Wang, {\em On the Smooth Compactification of
Siegel Spaces }, J. Diff. Geometry, {\bf 38}(1993), 351-386.

\bibitem{} W. Wang, {\em On the moduli space of principally
polarized abelian varieties}, Contemporary Math., {\bf 150}(1993),
361-365.

\bibitem{} H. Weber, {\em Theorie der Abel'schen Funktionen von
Geschlecht $3$}, Berlin, (1876).

\bibitem{} R. Weissauer, {\em Vektorwertige Siegelsche
Modulformen kleine Gewichtes}, J. reine angew. Math., {\bf 343}(1983),
184-202.

\bibitem{} R. Weissauer, {\em Untervariet{\"a}ten der Siegelschen
Modulmannigfatigkeiten von allgemeinen Typ}, Math. Ann., {\bf 343}(1983),
209-220.

\bibitem{} R. Weissauer, {\em Differentialformen zu Untergruppen der
Siegelschen Modulgruppe zweiten Grades}, J. reine angew Math., {\bf 391}(1988),
100-156.

\bibitem{} T. Yamazaki, {\em Jacobi forms and a Maass relation for
Eisenstein series}, J. Fac. Sci. Univ. Tokyo, {\bf 33}(1986),
295-310.

\bibitem{} J.-H. Yang, {\em Harmonic Analysis on the Quotient Spaces of
Heisenberg Groups}, Nagoya Math. J., {\bf 123}(1991), 103-117.

\bibitem{} J.-H. Yang, {\em Harmonic Analysis on the Quotient Spaces of
Heisenberg Groups II }, J. Number Theory, {\bf 49}(1)(1994), 63-72.

\bibitem{} J.-H. Yang ,{\em A decomposition theorem on
differential polynomials of theta functions of high level}, Japanese J. of
Mathematics, The Mathematical Society of Japan, New Series, {\bf 22}(1)(1996),
37-49.

\bibitem{} J.-H. Yang, {\em Lattice Representations of the Heisenberg Group
$H_{\BR}^{(g,h)}$}, Math. Annalen, {\bf 317}(2000), 309-323.

\bibitem{} J.-H. Yang, {\em Unitary Representations of the
Heisenberg Group $H_{\BR}^{(g,h)}$}, Proc. of Workshops in
Pure Mathematics on Groups and Algebraic Structures, the Korean
Academic Council, {\bf 5}, Part I, (1996), 144-159.

\bibitem{} J.-H. Yang, {\em Hermite Differential Operators
}, submited.

\bibitem{} J.-H. Yang, {\em Fock Representations}, J. Korean Math. Soc., {\bf 34}(1997),
345-370.

\bibitem{} J.-H. Yang, {\em On Theta Functions},
Kyungpook Math. J., {\bf 35 }(1996), 857-875.

\bibitem{} J.-H. Yang, {\em Geometrical theory of
Siegel modular forms}, Proc. Topology and Geometry Research
Center, {\bf 1}(1990), 73-83.

\bibitem{} J.-H. Yang, {\em Some Results on Jacobi
Forms of Higher Degree}, RIMS Symp. on Automorphic Forms
and Associated Zeta Functions, RIMS K${\hat o}$ky${\hat u}$-roku, {\bf 805}(1992), 36-50.

\bibitem{} J.-H. Yang, {\em Vanishing theorems on Jacobi
forms of higher degree}, J. Korean Math. Soc., {\bf 30}(1)(1993),
185-198.

\bibitem{} J.-H. Yang, {\em On Jacobi forms of higher degree (Korean)}, the Notices of the
Pyungsan Institute for Mathematical Sciences, {\bf 3}(1993), 3-29.

\bibitem{} J.-H. Yang, {\em Stable Jacobi forms}, Proc. of Workshops in Pure Mathematics on Number Theory and
Algebra, the Korean Academic Council ; {\bf 13}, Part I, (1993),
31-51.

\bibitem{} J.-H. Yang, {\em The Siegel-Jacobi Operator}, Abh. Math. Sem. Univ. Hamburg, {\bf 63}(1993), 135-146.

\bibitem{} J.-H. Yang, {\em Remarks on Jacobi forms of higher
degree}, Proc. of the 1993 Workshop on Automorphic Forms and
Related Topics, edited by Jin-Woo Son and Jae-Hyun Yang, the
Pyungsan Institute for Mathematical Sciences, (1993), 33-58.

\bibitem{} J.-H. Yang, {\em Singular Jacobi
Forms}, Trans. of American Math. Soc., {\bf 347}(6)(1995),
2041-2049.

\bibitem{} J.-H. Yang, {\em Construction of Modular Forms from Jacobi Forms}, Canadian J.
of Math., {\bf 47}(6)(1995), 1329-1339.

\bibitem{} J.-H. Yang, {\em Kac-Moody Algebras, the Monstrous
Moonshine, Jacobi Forms and Infinite Products}, Proceedings of
the 1995 Symposium on Number Theory, Geometry and Related Topics,
edited by Jin-Woo Son and Jae-Hyun Yang, the Pyungsan Institute
for Mathematical Sciences, (May 1996), 13-82.

\bibitem{} J.-H. Yang, {\em Stable Automorphic Forms},
Proceedings of Japan-Korea Joint Seminar on Transcendental Number
Theory and Related Topics, (1998), 101-126.

\bibitem{} J.-H. Yang, {\em Invariant metrics on the Siegel-Jacobi
space}, preprint (PIMS 99-4), the Pyungsan Institute for
Mathematical Sciences, (1999).

\bibitem{} J.-H. Yang, {\em A note on Maass-Jacobi forms
}, preprint (PIMS 2000-1), the Pyungsan Institute for
Mathematical Sciences, (2000).

\bibitem{} J.-H. Yang and Jin Woo Son, {\em A
note on Jacobi forms of higher degree}, J. Korean Math. Soc., {\bf 28}(2)(1991),
341-358.

\bibitem{} D. Zagier, {\em Note on the Landweber-Stong Elliptic Genus},
LNM, Springer-Verlag, {\bf 1326}(1988), 216-224.

\bibitem{} D. Zagier, {\em Periods of modular forms and Jacobi theta
functions}, Invent. Math., {\bf 104}(1991), 449-465.

\bibitem{} C. Ziegler, {\em Jacobi Forms of Higher Degree}, Abh. Math.
Sem. Univ. Hamburg, {\bf 59}(1989), 191-224.


\end{thebibliography}
\end{document}